\documentclass[11pt,reqno]{amsart}
\usepackage{amsmath,amssymb}
\usepackage[colorlinks=true,linkcolor=magenta,citecolor=blue]{hyperref}
\usepackage[margin=1in]{geometry}
\usepackage{enumitem} 
\usepackage{mathabx}
\usepackage{epsfig}

\usepackage{graphicx}
\usepackage{mathrsfs}
\usepackage{color}
\usepackage{makecell,multirow,diagbox}
\usepackage{amsfonts}
\usepackage{amsmath,amsxtra,latexsym,amsthm, amssymb, amscd, pb-diagram}

\newtheorem{dn}{Definition}[section]
\newtheorem{dl}{Theorem}[section]
\newtheorem{md}{Proposition}[section]
\newtheorem{bd}{Lemma}[section]
\newtheorem{hq}{Corollary}[section]
\newtheorem{nx}{Remark}[section]
\newtheorem{vd}{Example}[section]
\newcommand{\R}{\mathbb{R}}

\newcommand{\e}{\varepsilon}
\newcommand{\ity}{\infty}
\newcommand{\f}{\frac}

\newcommand{\bbd}{\begin{bd}}
\newcommand{\ebd}{\end{bd}}
\newcommand{\bdn}{\begin{dn}}
\newcommand{\edn}{\end{dn}}
\newcommand{\bhq}{\begin{hq}}
\newcommand{\ehq}{\end{hq}}
\newcommand{\bdl}{\begin{dl}}
\newcommand{\edl}{\end{dl}}
\newcommand{\bnx}{\begin{nx}}
\newcommand{\enx}{\end{nx}}
\newcommand{\bmd}{\begin{md}}
\newcommand{\emd}{\end{md}}
\newcommand{\bvd}{\begin{vd}}
\newcommand{\evd}{\end{vd}}

\usepackage[displaymath, mathlines]{lineno}

\title[Semi-linear damped $\sigma$-evolution equations with a modulus of continuity term]{A remark on semi-linear damped $\sigma$-evolution equations with a modulus of continuity term in nonlinearity}
\author{Tuan Anh Dao}
\address{Tuan Anh Dao \hfill\break
$\quad$ School of Applied Mathematics and Informatics, Hanoi University of Science and Technology, No.1 Dai Co Viet road, Hanoi, Vietnam \hfill\break
Faculty for Mathematics and Computer Science, TU Bergakademie Freiberg, Pr\"{u}ferstr. 9, 09596, Freiberg, Germany}
\email{anh.daotuan@hust.edu.vn}
\author{Michael Reissig}
\address{Michael Reissig \hfill\break
$\quad$ Faculty for Mathematics and Computer Science, TU Bergakademie Freiberg, Pr\"{u}ferstr. 9, 09596, Freiberg, Germany}
\email{reissig@math.tu-freiberg.de}

\begin{document}
	
\begin{abstract}
In this article, we indicate that under suitable assumptions of a modulus of continuity we obtain either the global (in time) existence of small data Sobolev solutions or the blow-up result of local (in time) Sobolev solutions to semi-linear damped $\sigma$-evolution equations with a modulus of continuity term in nonlinearity.
\end{abstract}

	\maketitle
	
	\tableofcontents
	
\section{Introduction and main results} \label{Sec.main}
In this paper we consider the following Cauchy problem for the semi-linear damped $\sigma$-evolution equations with modulus of continuity term in nonlinearity:
\begin{equation} \label{pt1.1}
\begin{cases}
u_{tt}+ (-\Delta)^\sigma u+ (-\Delta)^{\delta} u_t= |u|^{p_1^*(m,n)}\mu(|u|), \\
u(0,x)= u_0(x),\,\,\, u_t(0,x)=u_1(x),
\end{cases}
\end{equation}
and
\begin{equation} \label{pt1.2}
\begin{cases}
u_{tt}+ (-\Delta)^\sigma u+ (-\Delta)^{\delta} u_t= |u_t|^{p_2^*(m,n)}\mu(|u_t|), \\
u(0,x)= u_0(x),\,\,\, u_t(0,x)=u_1(x),
\end{cases}
\end{equation}
where $\sigma \ge 1$, $\delta\in [0,\frac{\sigma}{2}]$, given real numbers $p_1^*(m,n):=  1+\frac{2m\sigma}{n- 2m\delta}$ and $p_2^*(m,n):=  1+\frac{m\sigma}{n}$, with $m \in [1,2)$ and $n \ge 1$. Here the function $\mu= \mu(|u|)$ stands for some moduli of continuity. We are interested in studying two main equations in the present paper including $\sigma$-evolution equation with frictional damping $\delta= 0$ and that with structurally damping $\delta \in (0,\frac{\sigma}{2})$. \medskip

There are several recent papers (see, for instance, \cite{DabbiccoEbert,DuongKainaneReissig,DuongReissig,DaoReissig}) concerning two of the most common cases of power nonlinearity $|\partial_t^j u|^p$ with $j= 0,\,1$ to the semi-linear damped $\sigma$-evolution equations, that is, to the following Cauchy problems:
\begin{equation} \label{pt1.3}
\begin{cases}
u_{tt}+ (-\Delta)^\sigma u+ (-\Delta)^{\delta} u_t= |\partial_t^j u|^p \\
u(0,x)= u_0(x),\,\,\, u_t(0,x)=u_1(x),
\end{cases}
\end{equation}
with $\sigma \ge 1$, $\delta\in [0,\frac{\sigma}{2}]$ and a real number $p>1$. In particular, the authors in \cite{DuongKainaneReissig,DuongReissig} used $(L^1\cap L^2)- L^2$ estimates and $L^2- L^2$ estimates for the solutions to the corresponding linear equations with vanishing right-hand side, i.e. the mixing of additional $L^1$ regularity for the data on the basis of $L^2- L^2$ estimates to prove the global (in time) existence of small data solutions to (\ref{pt1.3}) in the cases $\delta= \frac{\sigma}{2}$ or $\delta= 0$. Meanwhile, a different strategy appearing in \cite{DabbiccoEbert} is to take account of additional $L^{\eta}\cap L^{\bar{q}}$ regularity, with small $\eta$ and large $\bar{q}$, in place of additional $L^1$ regularity. On the one hand, this strategy gives the global (in time) existence of small data solutions to the semi-linear models to (\ref{pt1.3}) in the case $\delta \in (0,\frac{\sigma}{2})$. Besides, some blow-up results were obtained in the latter paper to really find critical exponent for $p$. Here, critical exponent $p_{crit}=p_{crit}(n)$ means that for some range of admissible $p>p_{crit}$, the so-called supercritical cases, there exists a global (in time) Sobolev solution for small initial data from a suitable function space. Moreover, one may find suitable small data such that there exists no global (in time) Sobolev solution in the subcritical cases $1< p \le p_{crit}$. In other words, we have, in general, only local (in time) Sobolev solutions under this assumption for the exponent $p$. Also, in \cite{DabbiccoEbert} the sharpness of the critical exponents to (\ref{pt1.3}) are given by $p_{crit}= p_1^*(1,n)= 1+ \frac{2\sigma}{n- 2\delta}$ with $j=0$ and $p_{crit}= p_2^*(1,n)= 1+ \frac{\sigma}{n}$ with $j=1$. For this reason, we can expect to look for global (in time) results to (\ref{pt1.1}) and (\ref{pt1.2}) by assuming additional $L^m$ regularity for the data, with $m \in [1,2)$, in the whole supercritical ranges. \smallskip

\textit{Modulus of continuity} is a well-known notation to describe the regularity of a function with respect to desired variables (see more \cite{CicognaniLorenz,Jah,Lorenz}). Some linear Cauchy problems with low regular coefficients combined with modulus of continuity were considered in these references to study the uniqueness and the conditional stability. In the present paper, we investigate power nonlinearity linked to some moduli of continuity as another more complicated type of nonlinearity terms. On the one hand, the main difficulty appearing is to deal with estimating modulus of continuity terms in our proofs. Nevertheless, considering the nonlinearities combined with modulus of continuity to (\ref{pt1.1}) brings some benefits to find $p_{crit}$. This means the connection is understood as an important approach to describe the behavior of critical exponent (see latter, Remarks \ref{remark1.2} and \ref{remark2.2}). More recently, in the recent paper of the second author and his collaborators \cite{EbertGirardiReissig} the authors focused on studying the global (in time) existence and the blow-up result of the solutions as well to (\ref{pt1.1}) in the case of classical damped wave $\sigma= 1$ and $\delta= 0$ by additional $L^1$ regularity for the data. As motivated from the cited paper, in the present paper we want to develop their techniques to further extend some results in the more generalized cases. Hence, our goal is twofold. The first motivation of this paper is to derive the global (in time) existence of small data Sobolev solutions to (\ref{pt1.1}) for all $\delta\in [0,\frac{\sigma}{2}]$ under a suitable assumption of moduli of continuity $\mu$ and by using additional $L^m$ regularity, with $m \in [1,2)$, for the data. Moreover, under the same assumption we not only indicate the global (in time) existence of small data energy solutions to (\ref{pt1.1}) but also (\ref{pt1.2}) in special case $\delta= \frac{\sigma}{2}$. The second motivation of this paper is to varify the critical exponents $p_{crit}= p_1^*(1,n)$ and $p_{crit}= p_2^*(1,n)$, respectively, to (\ref{pt1.1}) and (\ref{pt1.2}) under a inverse assumption of moduli of continuity $\mu$ when $\sigma \ge 1$ and $\delta \in [0,\frac{\sigma}{2}]$ are integers. \medskip

$\qquad$ \textbf{Main results} \medskip

First we state the global (in time) existence of small data Sobolev solutions to (\ref{pt1.1}) in the case $\delta \in [0,\frac{\sigma}{2}]$.
\bdl \label{dlexistence1}
Let $\sigma \ge 1$, $\delta\in [0,\frac{\sigma}{2}]$ and $m \in [1,2)$. Let $0< r\le \sigma$. We assume the space dimensions satisfying the conditions
\begin{equation}
\begin{cases}
2m_0\delta< n< 2r \qquad\, \text{ if } \delta \in [0,\frac{\sigma}{2}), \\
m\sigma< n< 2r \qquad\quad \text{ if } \delta= \frac{\sigma}{2},
\end{cases}
\label{spacedimensionCon1}
\end{equation}
with $\frac{1}{m_0}=\frac{1}{m}- \frac{1}{2}$ and the following assumptions of modulus of continuity:
\begin{equation}
s\mu'(s) \lesssim \mu(s), \label{*}
\end{equation}
and
\begin{equation}
\int_{C_0}^\ity \frac{\mu(\frac{1}{s})}{s}ds < \ity, \label{existenceCon}
\end{equation}
with a sufficiently lare constant $C_0 >0$. Then, there exists a constant $\e>0$ such that for any small data
$$ (u_0,u_1) \in \big(L^m \cap H^r \big) \times \big(L^m \cap L^2\big) $$
satisfying the assumption
\[ \|u_0\|_{L^m \cap H^r}+ \|u_1\|_{L^m \cap L^2} \le \e, \]
we have a uniquely determined global (in time) small data Sobolev solution
$$ u \in C([0,\ity),H^r) $$
to (\ref{pt1.1}). Moreover, the following estimates hold:
\begin{align*}
\|u(t,\cdot)\|_{L^2}& \lesssim (1+t)^{-\frac{n}{2(\sigma -\delta)}(\frac{1}{m}-\frac{1}{2})+ \frac{\delta}{\sigma -\delta}} \big(\|u_0\|_{L^m \cap H^r}+ \|u_1\|_{L^m \cap L^2}\big), \\
\big\||D|^r u(t,\cdot)\big\|_{L^2}& \lesssim (1+t)^{-\frac{n}{2(\sigma -\delta)}(\frac{1}{m}-\frac{1}{2})- \frac{r- 2\delta}{2(\sigma -\delta)}} \big(\|u_0\|_{L^m \cap H^r}+ \|u_1\|_{L^m \cap L^2}\big).
\end{align*}
\edl

Next we state the global (in time) existence of small data energy solutions to (\ref{pt1.1}) in the case $\delta= \frac{\sigma}{2}$.

\bdl \label{dlexistence2}
Let $\sigma \ge 1$, $\delta= \frac{\sigma}{2}$ and $m \in [1,2)$. We assume the following condition of space dimensions:
\begin{equation}
m\sigma< n< 2\sigma,
\label{spacedimensionCon2}
\end{equation}
and modulus of continuity satisfies the assumptions (\ref{*}) and (\ref{existenceCon}) as in Theorem \ref{dlexistence1}. Then, there exists a constant $\e>0$ such that for any small data
$$ (u_0,u_1) \in \big(L^m \cap H^\sigma \big) \times \big(L^m \cap L^2\big) $$
satisfying the assumption
\[ \|u_0\|_{L^m \cap H^\sigma}+ \|u_1\|_{L^m \cap L^2} \le \e, \]
we have a uniquely determined global (in time) small data energy solution
$$ u \in C([0,\ity),H^\sigma) \cap C^1([0,\ity),L^2) $$
to (\ref{pt1.1}). Moreover, the following estimates hold:
\begin{align*}
\|u(t,\cdot)\|_{L^2}& \lesssim (1+t)^{-\frac{n}{\sigma}(\frac{1}{m}-\frac{1}{2})+ 1} \big(\|u_0\|_{L^m \cap H^\sigma}+ \|u_1\|_{L^m \cap L^2}\big), \\
\big\|\big(|D|^\sigma u(t,\cdot), u_t(t,\cdot)\big)\big\|_{L^2}& \lesssim (1+t)^{-\frac{n}{\sigma}(\frac{1}{m}-\frac{1}{2})} \big(\|u_0\|_{L^m \cap H^\sigma}+ \|u_1\|_{L^m \cap L^2}\big).
\end{align*}
\edl

Now we state the global (in time) existence of large regular solutions to (\ref{pt1.2}) in the case $\delta= \frac{\sigma}{2}$.

\bdl \label{dlexistence3}
Let $\sigma \ge 1$, $\delta= \frac{\sigma}{2}$ and $m \in [1,2)$. Let $r> \sigma+ \frac{n}{2}$. We assume the following condition:
\begin{equation}
r\le  2\sigma- \frac{n}{m_0},
\label{spacedimensionCon3}
\end{equation}
with $\frac{1}{m_0}=\frac{1}{m}- \frac{1}{2}$ and modulus of continuity satisfies the assumptions (\ref{*}) and (\ref{existenceCon}) as in Theorem \ref{dlexistence1}. Then, there exists a constant $\e>0$ such that for any small data
$$ (u_0,u_1) \in \big(L^m \cap H^r \big) \times \big(L^m \cap H^{r-\sigma}\big) $$
satisfying the assumption
\[ \|u_0\|_{L^m \cap H^r}+ \|u_1\|_{L^m \cap H^{r-\sigma}} \le \e, \]
we have a uniquely determined global (in time) small data energy solution
$$ u \in C([0,\ity),H^r) \cap C^1([0,\ity),H^{r-\sigma}) $$
to (\ref{pt1.2}). Moreover, the following estimates hold:
\begin{align*}
\|u(t,\cdot)\|_{L^2}& \lesssim (1+t)^{-\frac{n}{\sigma}(\frac{1}{m}-\frac{1}{2})+ 1} \big(\|u_0\|_{L^m \cap H^r}+ \|u_1\|_{L^m \cap H^{r-\sigma}}\big), \\
\|u_t(t,\cdot)\|_{L^2}& \lesssim (1+t)^{-\frac{n}{\sigma}(\frac{1}{m}-\frac{1}{2})} \big(\|u_0\|_{L^m \cap H^r}+ \|u_1\|_{L^m \cap H^{r-\sigma}}\big), \\
\big\|\big(|D|^r u(t,\cdot), |D|^{r- \sigma}u_t(t,\cdot)\big)\big\|_{L^2}& \lesssim (1+t)^{-\frac{n}{\sigma}(\frac{1}{m}-\frac{1}{2})- \frac{r- \sigma}{\sigma}} \big(\|u_0\|_{L^m \cap H^r}+ \|u_1\|_{L^m \cap H^{r-\sigma}}\big).
\end{align*}
\edl

Now we state the blow-up result to (\ref{pt1.1}) in the case $\delta \in [0,\frac{\sigma}{2}]$.
\begin{dl} \label{dloptimal1}
Let $\sigma \ge 1$ and $\delta \in [0,\frac{\sigma}{2}]$ be integer numbers. We assume that
the first datum $u_0= 0$, whereas the second datum $u_1 \in L^1 \cap L^2$ satisfies the following relation:
\begin{equation} \label{optimaldata}
\int_{\R^n} u_1(x)dx >0.
\end{equation}
Moreover, we suppose the following assumption of modulus of continuity
\begin{equation}
\int_{C_0}^\ity \frac{\mu(\frac{1}{s})}{s}ds= \ity, \label{blowupCon}
\end{equation}
with a sufficiently lare constant $C_0 >0$. Then, there is no global (in time) Sobolev solution to (\ref{pt1.1}).
\end{dl}

Finally, we state the blow-up result to (\ref{pt1.2}) in the case $\delta= \frac{\sigma}{2}$.
\begin{dl} \label{dloptimal2}
Let $\sigma \ge 1$ and $\delta= \frac{\sigma}{2}$ be integer numbers. We assume that the first datum $u_0= 0$, whereas the second datum $u_1 \in L^1 \cap L^2$ satisfies the following relation:
\begin{equation} \label{optimaldata2}
\int_{\R^n} u_1(x)dx >0.
\end{equation}
Moreover, we suppose the assumption of modulus of continuity (\ref{blowupCon}) is fulfilled as in Theorem \ref{dloptimal1}. Then, there is no global (in time) energy solution to (\ref{pt1.2}).
\end{dl}

\begin{nx}
\fontshape{n}
\selectfont
If we plug $m=1$ in Theorems \ref{dlexistence1}, \ref{dlexistence2} and \ref{dlexistence3} then from Theorems \ref{dloptimal1} and \ref{dloptimal2} it is clear that the critical exponents $p_{crit}$ to (\ref{pt1.1}) is given by $p_{crit}= p_1^*(1,n)= 1+ \frac{2\sigma}{n- 2\delta}$ and $p_{crit}$ to (\ref{pt1.2}) is given by $p_{crit}= p_2^*(1,n)= 1+ \frac{\sigma}{n}$, i.e our results are optimal.
\end{nx}

\bnx
\fontshape{n}
\selectfont
In comparison with the known result from the paper \cite{Takeda} concerning the power nonlinearity $|u|^p$ to (\ref{pt1.1}), we want to underline that the obtained critical exponent $p_{crit}$ from Theorem \ref{dloptimal1} coincides with that in the cited paper.
\enx

\begin{nx} \label{remark1.2}
\fontshape{n}
\selectfont
From the assumption (\ref{existenceCon}) in Theorems \ref{dlexistence1}, \ref{dlexistence2}, \ref{dlexistence3} and the assumption (\ref{blowupCon}) in Theorems \ref{dloptimal1}, \ref{dloptimal2} here we want to emphasize that modulus of continuity comes into play to guarantee either the global (in time) existence or the blow-up result under its suitable conditions.
\end{nx}

\noindent \textbf{The organization of this paper} is presented as follows: In Section \ref{Sec.Pre}, we collect some basic properties about modulus of continuity and the $(L^m \cap L^2)- L^2$ estimates and the $L^2- L^2$, with $m\in [1,2)$, for the solutions to the linear Cauchy problems as well. We shall apply these estimates to prove our main results in Section \ref{Sec.Proof}. In particular, we give the detail proofs of the global (in time) existence of small data Sobolev solutions and small data energy solutions to (\ref{pt1.1}) and (\ref{pt1.2}) in Section \ref{Sec.Existence}. Moreover, in Section \ref{Sec.Optimal} we indicate the optimality of the power exponents.

\section{Preliminaries} \label{Sec.Pre}
In this section, we collect some preliminary knowledge needed in our proofs.
\subsection{Modulus of continuity} \label{Sec.Modulus}
First of all, we recall the signifficant properties and some typical examples about modulus of continuity as well (see, for instance, \cite{CicognaniLorenz,Lorenz} and the references therein).
\bdn \label{DefModulsofContinuity}
\fontshape{n}
\selectfont
Let $\mu: [0,c] \to [0,c]$ be a continuous, concave an increasing function with a sufficiently small positive constant $c$. Then $\mu$ is called a modulus of continuity if it satisfies
$$ \mu(0)= 0. $$
\edn

\begin{nx}
\fontshape{n}
\selectfont
Here we want to underline that in some sense the above definition of a modulus of continuity could be extended into $[0,1]$ or $[0,\ity)$ instead of $[0,c]$. However, the essential point of any modulus of continuity is its behavior near $0$. This property also comes into play in comparison between two different moduli of continuity. For this reason, it is sufficient to only consider moduli of continuity on intervals $[0,c]$ with a sufficiently small positive constant $c$.
\end{nx}

\begin{vd}
\fontshape{n}
\selectfont
In the following list, we give some typical examples about moduli of continuity which are arranged according to their regularity from the highest one to the lowest one: \medskip

\begin{tabular}{ll}
\hline
modulus of continuity &$\qquad$ frequently called name  \\
\hline
$\mu(s)= s$ &$\qquad$ Lipschitz-continuity \\
\smallskip
$\mu(s)= s\Big(\log\Big(\frac{1}{s}\Big)+ 1\Big)$ &$\qquad$ Log-Lip-continuity \\
\smallskip
$\mu(s)= s\Big(\log\Big(\frac{1}{s}\Big)+ 1\Big)\log^{[m]}\Big(\frac{1}{s}\Big), \quad m \ge 1$ &$\qquad$ Log-Log$^{[m]}$-Lip-continuity \\
\smallskip
$\mu(s)= s^\alpha, \quad \alpha \in (0,1)$ &$\qquad$ H\"{o}lder-continuity \\
\smallskip
$\mu(s)= \Big(\log\Big(\frac{1}{s}\Big)+ 1\Big)^{-\alpha}, \quad \alpha \in (0,\ity)$ &$\qquad$ Log$^{-\alpha}$-continuity \\
\hline
\end{tabular}
\end{vd}
where $\log^{[m]}(x)= \log\big(\log^{[m-1]}(x)\big)+ 1$ for $m \ge 2$, and $\log^{[1]}(x)= \log(x)+ 1$.
\begin{nx} \label{remark2.2}
\fontshape{n}
\selectfont
We can easily check that the moduli of continuity such as Lipschitz-continuity, Log-Lip-continuity, Log-Log$^{[m]}$-Lip-continuity, H\"{o}lder-continuity and Log$^{-\alpha}$-continuity with $\alpha \in (1,\ity)$ fulfill the assumption (\ref{existenceCon}) in Theorem \ref{dlexistence1}, whereas Log$^{-\alpha}$-continuity with $\alpha \in (0,1]$ satisfies the assumption (\ref{blowupCon}) in Theorem \ref{dloptimal1}.
\end{nx}

\subsection{Linear estimates} \label{Sec.Linear}
Let us consider the corresponding linear models with vanishing right-hand side in the following form:
\begin{equation}
u_{tt}+ (-\Delta)^\sigma u+ (-\Delta)^{\delta} u_t=0,\,\, u(0,x)= u_0(x),\,\, u_t(0,x)= u_1(x), \label{pt2.1}
\end{equation}
with $\sigma \ge 1$ and $\delta\in [0,\frac{\sigma}{2}]$. Main goal of this section is to collect $(L^m \cap L^2)- L^2$ and $L^2- L^2$ estimates for the solutions and some of their derivatives to (\ref{pt1.2}). These estimates play a crucial role to prove the global (in time) existence result to (\ref{pt1.1}) in the next section.

\subsubsection{The structural damping case $\delta \in (0,\frac{\sigma}{2}]$} 
We obtained the following estimates from the previous paper \cite{Dao} of the first author.
\bmd[Proposition 2.1 in \cite{Dao}] \label{md2.1}
Let $\delta=\f{\sigma}{2}$ in (\ref{pt2.1}) and $m \in [1,2)$. The Sobolev solutions to (\ref{pt2.1}) satisfy the $(L^m \cap L^2)-L^2$ estimates
\begin{align*}
\big\|\partial_t^j |D|^a u(t,\cdot)\big\|_{L^2} &\lesssim (1+t)^{-\frac{n}{\sigma}(\frac{1}{m}-\frac{1}{2})- \frac{a}{\sigma}-j}\|u_0\|_{L^m \cap H^{a+j\sigma}} \\ 
&\qquad \qquad + (1+t)^{1-\frac{n}{\sigma}(\frac{1}{m}-\frac{1}{2})- \frac{a}{\sigma}-j}\|u_1\|_{L^m \cap H^{[a+(j-1)\sigma]^+}},
\end{align*}
and the $L^2-L^2$ estimates
$$ \big\|\partial_t^j |D|^a u(t,\cdot)\big\|_{L^2} \lesssim (1+t)^{- \frac{a}{\sigma}-j}\|u_0\|_{H^{a+j\sigma}}+ (1+t)^{1- \frac{a}{\sigma}-j}\|u_1\|_{H^{[a+(j-1)\sigma]^+}}, $$
for any $a\ge 0$, $j=0,1$ and for all space dimensions $n \ge 1$.
\emd

\bmd[Proposition 2.2 in \cite{Dao}] \label{md2.2}
Let $\delta \in (0,\f{\sigma}{2})$ in (\ref{pt2.1}) and $m \in [1,2)$. The Sobolev solutions to (\ref{pt2.1}) satisfy the $(L^m \cap L^2)-L^2$ estimates
\begin{align*}
\big\|\partial_t^j |D|^a u(t,\cdot)\big\|_{L^2} &\lesssim (1+t)^{-\frac{n}{2(\sigma-\delta)}(\frac{1}{m}-\frac{1}{2})- \frac{a+2j\delta}{2(\sigma-\delta)}}\|u_0\|_{L^m \cap H^{a+j\sigma}}\\ 
&\qquad \qquad + (1+t)^{1-\frac{n}{2(\sigma-\delta)}(\frac{1}{m}-\frac{1}{2})- \frac{a+2j\delta}{2(\sigma-\delta)}}\|u_1\|_{L^m \cap H^{[a+(j-1)\sigma]^+}},
\end{align*}
and the $L^2-L^2$ estimates
$$ \big\|\partial_t^j |D|^a u(t,\cdot)\big\|_{L^2} \lesssim (1+t)^{-\frac{a+2j\delta}{2(\sigma-\delta)}}\|u_0\|_{H^{a+j\sigma}}+ (1+t)^{1- \frac{a+2j\delta}{2(\sigma-\delta)}}\|u_1\|_{H^{[a+(j-1)\sigma]^+}}, $$
for any $a\ge 0$, $j=0,1$ and for all space dimensions $n \ge 1$.
\emd

\begin{nx}
\fontshape{n}
\selectfont
In Proposition \ref{md2.2} we state estimates for the solution and some its derivatives to (\ref{pt2.1}) which hold for any space dimensions $n \ge 1$. Moreover, we may prove a better result under a restriction to space dimensions $n>2m_0 \delta$ with $m_0$ satisfying $\frac{1}{m_0}=\frac{1}{m}- \frac{1}{2}$. We get the following sharper estimates.
\end{nx}

\bmd[Proposition 2.3 in \cite{Dao}] \label{md2.3}
Let $\delta \in (0,\f{\sigma}{2})$ in (\ref{pt2.1}) and $m \in [1,2)$. The Sobolev solutions to (\ref{pt2.1}) satisfy the $(L^m \cap L^2)-L^2$ estimates
\begin{align*}
\big\|\partial_t^j |D|^a u(t,\cdot)\big\|_{L^2} &\lesssim (1+t)^{-\frac{n}{2(\sigma-\delta)}(\frac{1}{m}-\frac{1}{2})- \frac{a}{2(\sigma-\delta)}-j}\|u_0\|_{L^m \cap H^{a+j\sigma}}\\ 
&\qquad \qquad + (1+t)^{-\frac{n}{2(\sigma-\delta)}(\frac{1}{m}-\frac{1}{2})- \frac{a-2\delta}{2(\sigma-\delta)}-j}\|u_1\|_{L^m \cap H^{[a+(j-1)\sigma]^+}},
\end{align*}
and the $L^2-L^2$ estimates
$$ \big\|\partial_t^j |D|^a u(t,\cdot)\big\|_{L^2} \lesssim (1+t)^{-\frac{a}{2(\sigma-\delta)}-j}\|u_0\|_{H^{a+j\sigma}}+ (1+t)^{- \frac{a-2\delta}{2(\sigma-\delta)}-j}\|u_1\|_{H^{[a+(j-1)\sigma]^+}}, $$
for any $a\ge 0$, $j=0,1$ and for all space dimensions $n>2m_0 \delta$.
\emd

\subsubsection{The frictional damping case $\delta=0$}
Our approach is based on the paper \cite{DuongReissig}. According to the treatment of Proposition $2.1$ in \cite{DuongReissig}, with minor modifications in the steps of the proofs we derive the following estimates.
\bmd \label{md2.4}
Let $\delta= 0$ in (\ref{pt2.1}) and $m \in [1,2)$. The Sobolev solutions to (\ref{pt2.1}) satisfy the $(L^m \cap L^2)-L^2$ estimates
\begin{align*}
\big\|\partial_t^j |D|^a u(t,\cdot)\big\|_{L^2} &\lesssim (1+t)^{-\frac{n}{2\sigma}(\frac{1}{m}-\frac{1}{2})-\frac{a}{2\sigma}-j}\|u_0\|_{L^m \cap H^{a+j\sigma}} \\
&\qquad \qquad + (1+t)^{-\frac{n}{2\sigma}(\frac{1}{m}-\frac{1}{2})-\frac{a}{2\sigma}-j}\|u_1\|_{L^m \cap H^{[a+(j-1)\sigma]^+}},
\end{align*}
and the $L^2-L^2$ estimates
$$ \big\|\partial_t^j |D|^a u(t,\cdot)\big\|_{L^2} \lesssim (1+t)^{-\frac{a}{2\sigma}-j} \|u_0\|_{H^{a+j\sigma}}+ (1+t)^{-\frac{a}{2\sigma}-j}\|u_1\|_{H^{[a+(j-1)\sigma]^+}}, $$
for any $a\ge 0$, $j=0,1$ and for all space dimensions $n \ge 1$.
\emd
Finally, combining Propositions \ref{md2.1}, \ref{md2.3} and \ref{md2.4} we may conclude the following estimates.
\bmd \label{md2.5}
Let $\delta \in [0,\f{\sigma}{2}]$ in (\ref{pt2.1}) and $m \in [1,2)$. The Sobolev solutions to (\ref{pt2.1}) satisfy the $(L^m \cap L^2)-L^2$ estimates
\begin{align*}
\big\|\partial_t^j |D|^a u(t,\cdot)\big\|_{L^2} &\lesssim (1+t)^{-\frac{n}{2(\sigma-\delta)}(\frac{1}{m}-\frac{1}{2})- \frac{a}{2(\sigma-\delta)}-j}\|u_0\|_{L^m \cap H^{a+j\sigma}}\\ 
&\qquad \qquad + (1+t)^{-\frac{n}{2(\sigma-\delta)}(\frac{1}{m}-\frac{1}{2})- \frac{a-2\delta}{2(\sigma-\delta)}-j}\|u_1\|_{L^m \cap H^{[a+(j-1)\sigma]^+}},
\end{align*}
and the $L^2-L^2$ estimates
$$ \big\|\partial_t^j |D|^a u(t,\cdot)\big\|_{L^2} \lesssim (1+t)^{-\frac{a}{2(\sigma-\delta)}-j}\|u_0\|_{H^{a+j\sigma}}+ (1+t)^{- \frac{a-2\delta}{2(\sigma-\delta)}-j}\|u_1\|_{H^{[a+(j-1)\sigma]^+}}, $$
for any $a\ge 0$ and $j=0,1$. Here we assume the following conditions of space dimensions:
$$ \begin{cases}
n>2m_0 \delta \quad\,\, \text{ if } \delta \in [0, \frac{\sigma}{2}), \\
n \ge 1 \qquad \quad \text{ if } \delta= \frac{\sigma}{2}.
\end{cases}$$
\emd

\begin{nx}
\fontshape{n}
\selectfont
Here we want to underline that the decay estimates for the solutions to (\ref{pt2.1}) from Proposition \ref{md2.3} are better than those from Proposition \ref{md2.2} in the case $\delta \in (0,\frac{\sigma}{2})$. Hence, it is reasonable to apply all statements from Proposition \ref{md2.5} in the steps of the proofs to our global (in time) existence results.
\end{nx}

\section{Proofs for main results} \label{Sec.Proof}
The ideas of the following proofs are based on the recent paper of the second author and his collaborators \cite{EbertGirardiReissig} in which the authors focused on their considerations to (\ref{pt1.1}) with $\sigma= 1$ and $\delta= 0$.

\subsection{Global existence results} \label{Sec.Existence}  
For simplicity, in the following proofs we use the abbreviations $p^*:= p_1^*(m,n)= 1+ \frac{2m\sigma}{n- 2m\delta}$ or $p^*:= p_2^*(m,n)= 1+ \frac{m\sigma}{n}$ to (\ref{pt1.1}) or (\ref{pt1.2}), respectively.
\subsubsection{Philosophy of our approach}
We choose the data space
$$(u_0,u_1) \in \mathcal{A}^{r}_m:= \big(L^m \cap H^r\big) \times \big(L^m \cap H^{[r-\sigma]^+}\big) $$
and introduce the solution space $X(t)$ with the norm
\begin{align*}
\|u\|_{X(t)}:= &\sup_{0\le \tau \le t} \Big( f_1(\tau)^{-1}\|u(\tau,\cdot)\|_{L^2} + f_2(\tau)^{-1}\big\||D|^r u(\tau,\cdot)\big\|_{L^2} \\ 
&\qquad \quad + g_1(\tau)^{-1}\|u_t(\tau,\cdot)\|_{L^2} + g_2(\tau)^{-1}\big\||D|^{r-\sigma} u_t(\tau,\cdot)\big\|_{L^2} \Big),
\end{align*}
where the weights $f_j(\tau)$ and $g_j(\tau)$, with $j=1,\,2$, are suitably choosen from decay estimates of the solutions of the corresponding linear Cauchy problems (\ref{pt1.2}). Denoting $K_0(t,x)$ and $K_1(t,x)$ as the fundamental solutions to (\ref{pt1.2}) we may write the solutions to (\ref{pt1.2}) in the following form:
$$ u^{ln}(t,x)=K_0(t,x) \ast_{x} u_0(x)+ K_1(t,x) \ast_{x} u_1(x). $$
Applying Duhamel's principle gives the formal implicit representation of the solutions to (\ref{pt1.1}) and (\ref{pt1.2}) as follows:
$$ u(t,x)= u^{ln}(t,x) + \int_0^t K_1(t-\tau,x) \ast_x f(|u(t,x)|) d\tau=: u^{ln}(t,x)+ u^{nl}(t,x), $$
where $f(|u(t,x)|)= |u(t,x)|^{p^*}\mu(|u(t,x)|)$ and $f(|u(t,x)|)= |u_t(t,x)|^{p^*}\mu(|u_t(t,x)|)$ to (\ref{pt1.1}) and (\ref{pt1.2}), respectively. We define for all $t>0$ the operator $N: \quad u \in X(t) \longrightarrow Nu \in X(t)$ by the formula
$$Nu(t,x)= u^{ln}(t,x)+ u^{nl}(t,x). $$
We will prove that the operator $N$ satisfies the following two inequalities:
\begin{align}
\|Nu\|_{X(t)} &\lesssim \|(u_0,u_1)\|_{\mathcal{A}^{r}_m}+ \|u\|^{p^*}_{X(t)}, \label{pt3.3}\\
\|Nu-Nv\|_{X(t)} &\lesssim \|u- v\|_{X(t)} \big(\|u\|^{p^*-1}_{X(t)}+ \|v\|^{p^*-1}_{X(t)}\big). \label{pt3.4}
\end{align}
Here we note that from the definition of the norm in $X(t)$, by replacing $a= 0$ and $a= r$ in the statements from Proposition \ref{md2.5} we may conclude
$$ \|u^{ln}\|_{X(t)} \lesssim \|(u_0,u_1)\|_{\mathcal{A}^{r}_{m}}. $$
For this reason, in order to complete the proof of (\ref{pt3.3}) we prove only the following inequality:
\begin{equation}
\|u^{nl}\|_{X(t)} \lesssim \|u\|^{p^*}_{X(t)}. \label{pt3.31}
\end{equation}
Then, applying Banach's fixed point theorem we obtain local (in time) existence results of large data solutions and global (in time) existence results of small data solutions as well.

\subsubsection{Proof of Theorem \ref{dlexistence1}.}
We choose the solution space
$$ X(t):= C([0,\ity), H^r), $$
and the following weights:
$$f_1(\tau)= (1+\tau)^{-\frac{n}{2(\sigma -\delta)}(\frac{1}{m}-\frac{1}{2})+ \frac{\delta}{\sigma -\delta}},\quad f_2(\tau)=(1+\tau)^{-\frac{n}{2(\sigma -\delta)}(\frac{1}{m}-\frac{1}{2})- \frac{r- 2\delta}{2(\sigma -\delta)}}, $$
and $g_1(\tau)= g_2(\tau) \equiv 0$. First, let us prove the inequality (\ref{pt3.31}). To control some estimates for $u^{nl}$, our strategy is to use the $(L^m \cap L^2)- L^2$ estimates from Proposition \ref{md2.5} to get the following estimates for $k=0,1$:
\begin{equation}
\big\||D|^{kr} u^{nl}(t,\cdot)\big\|_{L^2} \lesssim \int_0^t (1+t-\tau)^{-\frac{n}{2(\sigma -\delta)}(\frac{1}{m}-\frac{1}{2})- \frac{kr- 2\delta}{2(\sigma -\delta)}}\big\||u(\tau,\cdot)|^{p^*}\mu\big(|u(\tau,\cdot)|\big)\big\|_{L^m \cap L^2}d\tau.
\label{t1.1}
\end{equation}
In order to estimate for $|u(\tau,\cdot)|^{p^*}\mu\big(|u(\tau,\cdot)|\big)$ in $L^m \cap L^2$, we proceed as follows:
\begin{equation}
\big\||u(\tau,\cdot)|^{p^*}\mu(|u(\tau,\cdot)|)\big\|_{L^m \cap L^2} \lesssim \big\||u(\tau,\cdot)|^{p^*}\big\|_{L^m \cap L^2} \big\|\mu\big(|u(\tau,\cdot)|\big)\big\|_{L^\ity}.
\label{t1.2}
\end{equation}
Since $\mu$ is an increasing function, we obtain
\begin{align*}
\big\|\mu\big(|u(\tau,\cdot)|\big)\big\|_{L^\ity} &\le \mu\big(\|u(\tau,\cdot)\|_{L^\ity}\big) \le \mu\big(C_1\|u(\tau,\cdot)\|_{\dot{H}^{r^*}}+ C_1\|u(\tau,\cdot)\|_{\dot{H}^r}\big) \\
&\le \mu\Big(C_1 (1+ \tau)^{-\frac{n}{2(\sigma -\delta)}(\frac{1}{m}-\frac{1}{2})- \frac{r^*- 2\delta}{2(\sigma -\delta)}} \|u\|_{X(\tau)}\Big),
\end{align*}
with a suitable positive constant $C_1$. Here we used the fractional Sobolev embedding from Proposition \ref{ProEmbedding} with the condition $r^*< \frac{n}{2}< r$. By choosing $r^*= \frac{n}{2}-\e$ with a sufficiently small positive number $\e$ and using the condition (\ref{spacedimensionCon1}) we get
\begin{equation}
\big\|\mu\big(|u(\tau,\cdot)|\big)\big\|_{L^\ity}\le \mu\Big(C_1 \e_0 (1+ \tau)^{-\frac{n}{2m(\sigma -\delta)}+ \frac{\delta}{\sigma -\delta}}\Big)= \mu\big(C_1 \e_0 (1+ \tau)^{-\alpha}\big),
\label{t1.3}
\end{equation}
where $\alpha:= \frac{n}{2m(\sigma -\delta)}- \frac{\delta}{\sigma -\delta}> 0$. For the other interesting term of (\ref{t1.2}), we re-write
$$ \big\||u(\tau,\cdot)|^{p^*}\big\|_{L^m \cap L^2}= \big\||u(\tau,\cdot)|^{p^*}\big\|_{L^m}+ \big\||u(\tau,\cdot)|^{p^*}\big\|_{L^2}= \|u(\tau,\cdot)\|^{p^*}_{L^{mp^*}}+ \|u(\tau,\cdot)\|^{p^*}_{L^{2p^*}}. $$
Employing the fractional Gagliardo-Nirenberg inequality from Proposition \ref{fractionalGagliardoNirenberg} gives
\begin{align*}
\|u(\tau,\cdot)\|^{p^*}_{L^{mp^*}} &\lesssim (1+\tau)^{-\frac{n}{2m(\sigma -\delta)}(p^*- 1)+ \frac{p^*\delta}{\sigma -\delta}}\|u\|^{p^*}_{X(\tau)},\\
\|u(\tau,\cdot)\|^{p^*}_{L^{2p^*}} &\lesssim (1+\tau)^{-\frac{np^*}{2(\sigma -\delta)}(\frac{1}{m}-\frac{1}{2p^*})+ \frac{p^* \delta}{\sigma -\delta}}\|u\|^{p^*}_{X(\tau)}.
\end{align*}
As a result, we derive
\begin{equation}
\big\||u(\tau,\cdot)|^{p^*}\big\|_{L^m \cap L^2} \lesssim (1+\tau)^{-\frac{n}{2m(\sigma -\delta)}(p^*- 1)+ \frac{p^*\delta}{\sigma -\delta}}\|u\|^{p^*}_{X(\tau)}= (1+\tau)^{-1}\|u\|^{p^*}_{X(\tau)},
\label{t1.4}
\end{equation}
where we note that $p^*= 1+ \frac{2m\sigma}{n- 2m\delta}$. Combining all the estimates from (\ref{t1.1}) to (\ref{t1.4}) we arrive at
\begin{equation}
\big\||D|^{kr} u^{nl}(t,\cdot)\big\|_{L^2} \lesssim \|u\|^{p^*}_{X(t)}\int_0^t (1+t-\tau)^{-\frac{n}{2(\sigma -\delta)}(\frac{1}{m}-\frac{1}{2})- \frac{kr- 2\delta}{2(\sigma -\delta)}}(1+\tau)^{-1} \mu\big(C_1 \e_0 (1+ \tau)^{-\alpha}\big)d\tau.
\label{t1.5}
\end{equation}
By splitting the above integral into two parts, on the one hand we derive the following estimate:
\begin{align}
I_1 &:= \int_0^{t/2}(1+t-\tau)^{-\frac{n}{2(\sigma -\delta)}(\frac{1}{m}-\frac{1}{2})- \frac{kr- 2\delta}{2(\sigma -\delta)}}(1+\tau)^{-1} \mu\big(C_1 \e_0 (1+ \tau)^{-\alpha}\big)d\tau \nonumber \\
&\,\,\, \lesssim (1+t)^{-\frac{n}{2(\sigma -\delta)}(\frac{1}{m}-\frac{1}{2})- \frac{kr- 2\delta}{2(\sigma -\delta)}}\int_0^{t/2}(1+\tau)^{-1} \mu\big(C_1 \e_0 (1+ \tau)^{-\alpha}\big)d\tau, \label{t1.6}
\end{align}
where we used $(1+t-\tau) \approx (1+t)$ for any $\tau \in [0,t/2]$. On the other hand, we can estimate the remaining integral as follows:
\begin{align*}
I_2 &:= \int_{t/2}^t (1+t-\tau)^{-\frac{n}{2(\sigma -\delta)}(\frac{1}{m}-\frac{1}{2})- \frac{kr- 2\delta}{2(\sigma -\delta)}}(1+\tau)^{-1} \mu\big(C_1 \e_0 (1+ \tau)^{-\alpha}\big)d\tau \\
&\,\,\, \lesssim (1+t)^{-\frac{n}{2(\sigma -\delta)}(\frac{1}{m}-\frac{1}{2})- \frac{kr- 2\delta}{2(\sigma -\delta)}}\int_{t/2}^t (1+t-\tau)^{-\frac{n}{2(\sigma -\delta)}(\frac{1}{m}-\frac{1}{2})- \frac{kr- 2\delta}{2(\sigma -\delta)}} \\
&\qquad \qquad \qquad \times (1+\tau)^{-1+ \frac{n}{2(\sigma -\delta)}(\frac{1}{m}-\frac{1}{2})+ \frac{kr- 2\delta}{2(\sigma -\delta)}}\mu\big(C_1 \e_0 (1+ \tau)^{-\alpha}\big)d\tau,
\end{align*}
where we used $(1+\tau) \approx (1+t) \text{ for any }\tau \in [t/2,t]$. Moreover, we pay attention that the following relation holds:
$$ -1+ \frac{n}{2(\sigma -\delta)}\Big(\frac{1}{m}-\frac{1}{2}\Big)+ \frac{kr- 2\delta}{2(\sigma -\delta)}< 0, $$
due to the condition (\ref{spacedimensionCon1}). It is clear to see that $1+\tau \ge 1+t- \tau$ for any $\tau \in [t/2,t]$. Consequently, we have
\begin{align*}
(1+ \tau)^{-\alpha} &\le (1+t- \tau)^{-\alpha}, \\
(1+\tau)^{-1+ \frac{n}{2(\sigma -\delta)}(\frac{1}{m}-\frac{1}{2})+ \frac{kr- 2\delta}{2(\sigma -\delta)}} &\le (1+t- \tau)^{-1+ \frac{n}{2(\sigma -\delta)}(\frac{1}{m}-\frac{1}{2})+ \frac{kr- 2\delta}{2(\sigma -\delta)}}.
\end{align*}
Hence, we arrive at
$$ I_2 \lesssim (1+t)^{-\frac{n}{2(\sigma -\delta)}(\frac{1}{m}-\frac{1}{2})- \frac{kr- 2\delta}{2(\sigma -\delta)}}\int_{t/2}^t (1+t-\tau)^{-1} \mu\big(C_1 \e_0 (1+t- \tau)^{-\alpha}\big)d\tau. $$
A standard change of variables leads to
\begin{equation}
I_2 \lesssim (1+t)^{-\frac{n}{2(\sigma -\delta)}(\frac{1}{m}-\frac{1}{2})- \frac{kr- 2\delta}{2(\sigma -\delta)}}\int_0^{t/2} (1+\tau)^{-1} \mu\big(C_1 \e_0 (1+ \tau)^{-\alpha}\big)d\tau.
\label{t1.7}
\end{equation}
From (\ref{t1.5}) to (\ref{t1.7}) we may conclude
\begin{align}
\big\||D|^{kr} u^{nl}(t,\cdot)\big\|_{L^2} &\lesssim (1+t)^{-\frac{n}{2(\sigma -\delta)}(\frac{1}{m}-\frac{1}{2})- \frac{kr- 2\delta}{2(\sigma -\delta)}}\|u\|^{p^*}_{X(t)} \int_0^{t/2} (1+\tau)^{-1} \mu\big(C_1 \e_0 (1+ \tau)^{-\alpha}\big)d\tau \nonumber \\
&\lesssim (1+t)^{-\frac{n}{2(\sigma -\delta)}(\frac{1}{m}-\frac{1}{2})- \frac{kr- 2\delta}{2(\sigma -\delta)}}\|u\|^{p^*}_{X(t)} \int_0^\ity (1+\tau)^{-1} \mu\big(C_1 \e_0 (1+ \tau)^{-\alpha}\big)d\tau. \label{t1.8}
\end{align}
Using again change of variable by denoting $s^{-1}:= C_1 \e_0 (1+ \tau)^{-\alpha}$ and a straightforward calculation give
\begin{equation}
\int_0^\ity (1+\tau)^{-1} \mu\Big(C_1 \e_0 (1+ \tau)^{-\alpha}\Big)d\tau \lesssim \int_{C_0}^\ity \frac{\mu\big(\frac{1}{s}\big)}{s}ds \lesssim 1,
\label{t1.9}
\end{equation}
where $C_0= \frac{1}{C_1 \e_0}$ is a sufficiently large constant. We notice that the assumption (\ref{existenceCon}) comes into play to guarantee the boundedness of the above integral. Therefore, from (\ref{t1.8}) and (\ref{t1.9}) we have shown that
$$ \big\||D|^{kr} u^{nl}(t,\cdot)\big\|_{L^2} \lesssim (1+t)^{-\frac{n}{2(\sigma -\delta)}(\frac{1}{m}-\frac{1}{2})- \frac{kr- 2\delta}{2(\sigma -\delta)}}\|u\|^{p^*}_{X(t)}. $$
From the definition of the norm in $X(t)$, we may conclude immediately the inequality (\ref{pt3.31}). \medskip

\noindent Next, let us prove the inequality (\ref{pt3.4}). For two elements $u$ and $v$ from $X(t)$, using again the $(L^m \cap L^2)- L^2$ estimates from Proposition \ref{md2.5} we get the following estimates with $k= 0,1$:
\begin{align*}
&\big\||D|^{kr}(Nu- Nv)(t,\cdot)\big\|_{L^2}= \big\||D|^{kr}(u^{nl}- v^{nl})(t,\cdot)\big\|_{L^2} \\
&\qquad \lesssim \int_0^t (1+t-\tau)^{-\frac{n}{2(\sigma -\delta)}(\frac{1}{m}-\frac{1}{2})- \frac{kr- 2\delta}{2(\sigma -\delta)}}\big\||u(\tau,\cdot)|^{p^*}\mu\big(|u(\tau,\cdot)|\big)- |v(\tau,\cdot)|^{p^*}\mu\big(|v(\tau,\cdot)|\big)\big\|_{L^m \cap L^2}d\tau.
\end{align*}
By using the mean value theorem we get the integral representation
\begin{align*}
&|u(\tau,x)|^{p^*}\mu\big(|u(\tau,x)|\big)- |v(\tau,x)|^{p^*}\mu\big(|v(\tau,x)|\big)\\ 
&\qquad= \big(u(\tau,x)- v(\tau,x)\big)\int_0^1 d_{|u|}H\big(\omega u(\tau,x)+(1-\omega)v(\tau,x)\big)\,d\omega,
\end{align*}
where $H(u)= |u|^{p^*}\mu(|u|)$. Thanks to the condition (\ref{*}) of modulus of continuity, we can estimate
$$ d_{|u|}H(u)= p^*|u|^{p^*- 1}\mu(|u|)+ |u|^{p^*}d_{|u|}\mu(|u|) \lesssim |u|^{p^*- 1}\mu(|u|). $$
Hence, we obtain
\begin{align*}
&|u(\tau,x)|^{p^*}\mu\big(|u(\tau,x)|\big)- |v(\tau,x)|^{p^*}\mu\big(|v(\tau,x)|\big)\\ 
&\qquad \lesssim \big(u(\tau,x)- v(\tau,x)\big)\int_0^1 \big|\omega u(\tau,x)+(1-\omega)v(\tau,x)\big|^{p^*- 1}\mu\big(\big|\omega u(\tau,x)+(1-\omega)v(\tau,x)\big|\big)\,d\omega.
\end{align*}
For this reason, we arrive at the following estimate:
\begin{align*}
&\big\||u(\tau,\cdot)|^{p^*}\mu\big(|u(\tau,\cdot)|\big)- |v(\tau,\cdot)|^{p^*}\mu\big(|v(\tau,\cdot)|\big)\big\|_{L^m \cap L^2}\\ 
&\,\,\, \lesssim \int_0^1 \big\|\big(u(\tau,x)- v(\tau,x)\big) \big|\omega u(\tau,x)+(1-\omega)v(\tau,x)\big|^{p^*- 1}\mu\big(\big|\omega u(\tau,x)+(1-\omega)v(\tau,x)\big|\big)\big\|_{L^m \cap L^2}\,d\omega.
\end{align*}
In the similar approach to the proof of the inequality (\ref{pt3.31}), we may conclude the inequality (\ref{pt3.4}). Summarizing, the proof is completed.

\subsubsection{Proof of Theorem \ref{dlexistence2}.}
We choose the solution space
$$ X(t):= C([0,\ity), H^\sigma) \cap C^1([0,\ity), L^2), $$
and the following weights:
$$f_1(\tau)= (1+\tau)^{-\frac{n}{\sigma}(\frac{1}{m}-\frac{1}{2})+ 1},\quad f_2(\tau)= g_1(\tau)= (1+\tau)^{-\frac{n}{\sigma}(\frac{1}{m}-\frac{1}{2})}, $$
and $g_2(\tau) \equiv 0$. First, let us prove the inequality (\ref{pt3.31}). In the same way as we did in the proof of Theorem \ref{dlexistence1}, we obtain the following estimates:
\begin{align}
\|u^{nl}(t,\cdot)\|_{L^2} &\lesssim (1+t)^{-\frac{n}{\sigma}(\frac{1}{m}-\frac{1}{2})+ 1} \|u\|^{p^*}_{X(t)}, \nonumber \\
\big\||D|^\sigma u^{nl}(t,\cdot)\big\|_{L^2} &\lesssim  (1+t)^{-\frac{n}{\sigma}(\frac{1}{m}-\frac{1}{2})} \|u\|^{p^*}_{X(t)}, \label{t2.1}
\end{align}
provided that condition (\ref{spacedimensionCon2}) is fulfilled. Analogously, controlling $\|u_t^{nl}(t,\cdot)\|_{L^2}$ will be handled as we did to get (\ref{t2.1}). Hence, we arrive at
$$ \|u_t^{nl}(t,\cdot)\|_{L^2} \lesssim  (1+t)^{-\frac{n}{\sigma}(\frac{1}{m}-\frac{1}{2})} \|u\|^{p^*}_{X(t)}. $$
From the definition of the norm in $X(t)$, we may conclude immediately the inequality (\ref{pt3.31}). \medskip

Next, let us prove the inequality (\ref{pt3.4}). Then, repeating the proof of Theorem \ref{dlexistence1} and using the same treatment as in the above steps, we may also conclude the inequality (\ref{pt3.4}). Summarizing, the proof is completed.

\subsubsection{Proof of Theorem \ref{dlexistence3}.}
We choose the solution space
$$ X(t):= C([0,\ity), H^r) \cap C^1([0,\ity), H^{r- \sigma}), $$
and the following weights:
$$f_1(\tau)= (1+\tau)^{-\frac{n}{\sigma}(\frac{1}{m}-\frac{1}{2})+ 1},\quad g_1(\tau)= (1+\tau)^{-\frac{n}{\sigma}(\frac{1}{m}-\frac{1}{2})},\quad f_2(\tau)= g_2(\tau)= (1+\tau)^{-\frac{n}{\sigma}(\frac{1}{m}-\frac{1}{2})-\frac{r-\sigma}{\sigma}}. $$
First, let us prove the inequality (\ref{pt3.31}). To control $u^{nl}$, we use the $(L^m \cap L^2)- L^2$ estimates from Proposition \ref{md2.5} to get the following estimates:
\begin{equation}
\|u^{nl}(t,\cdot)\|_{L^2} \lesssim \int_0^t (1+t-\tau)^{-\frac{n}{\sigma}(\frac{1}{m}-\frac{1}{2})+1}\big\||u_t(\tau,\cdot)|^{p^*}\mu\big(|u_t(\tau,\cdot)|\big)\big\|_{L^m \cap L^2}d\tau.
\label{t3.1}
\end{equation}
We derive
\begin{equation}
\big\||u_t(\tau,\cdot)|^{p^*}\mu(|u_t(\tau,\cdot)|)\big\|_{L^m \cap L^2} \lesssim \big\||u_t(\tau,\cdot)|^{p^*}\big\|_{L^m \cap L^2} \big\|\mu\big(|u_t(\tau,\cdot)|\big)\big\|_{L^\ity}.
\label{t3.2}
\end{equation}
Because $\mu$ is an increasing function, we obtain
\begin{align}
\big\|\mu\big(|u_t(\tau,\cdot)|\big)\big\|_{L^\ity} &\le \mu\big(\|u_t(\tau,\cdot)\|_{L^\ity}\big) \le \mu\big(C_1\|u_t(\tau,\cdot)\|_{\dot{H}^{r^*}}+ C_1\|u_t(\tau,\cdot)\|_{\dot{H}^{r-\sigma}}\big) \nonumber \\
&\le \mu\Big(C_1 (1+ \tau)^{-\frac{n}{\sigma}(\frac{1}{m}-\frac{1}{2})- \frac{r^*}{\sigma}} \|u\|_{X(\tau)}\Big) \le \mu\big(C_1 \e_0 (1+ \tau)^{-\alpha}\big), \label{t3.3}
\end{align}
with a suitable positive constant $C_1$ and $\alpha:= \frac{n}{\sigma}(\frac{1}{m}-\frac{1}{2})> 0$. Here we used the fractional Sobolev embedding from Proposition \ref{ProEmbedding} with the condition $r^*< \frac{n}{2}< r- \sigma$. For the other interesting term of (\ref{t1.2}), we re-write
$$ \big\||u_t(\tau,\cdot)|^{p^*}\big\|_{L^m \cap L^2}= \big\||u_t(\tau,\cdot)|^{p^*}\big\|_{L^m}+ \big\||u_t(\tau,\cdot)|^{p^*}\big\|_{L^2}= \|u_t(\tau,\cdot)\|^{p^*}_{L^{mp^*}}+ \|u_t(\tau,\cdot)\|^{p^*}_{L^{2p^*}}. $$
Applying the fractional Gagliardo-Nirenberg inequality from Proposition \ref{fractionalGagliardoNirenberg} we obtain
\begin{equation}
\big\||u_t(\tau,\cdot)|^{p^*}\big\|_{L^m \cap L^2} \lesssim (1+\tau)^{-\frac{n}{m\sigma}(p^*- 1)}\|u\|^{p^*}_{X(\tau)}= (1+\tau)^{-1}\|u\|^{p^*}_{X(\tau)},
\label{t3.4}
\end{equation}
where we note that $p^*= 1+ \frac{m\sigma}{n}$. Combining all the estimates from (\ref{t3.1}) to (\ref{t3.4}) we arrive at
\begin{equation}
\|u^{nl}(t,\cdot)\|_{L^2} \lesssim \|u\|^{p^*}_{X(t)}\int_0^t (1+t-\tau)^{-\frac{n}{\sigma}(\frac{1}{m}-\frac{1}{2})+1}(1+\tau)^{-1} \mu\big(C_1 \e_0 (1+ \tau)^{-\alpha}\big)d\tau.
\label{t3.5}
\end{equation}
By the similar arguments as in the proof of Theorem \ref{dlexistence1} we may conclude the following estimate:
\begin{equation}
\|u^{nl}(t,\cdot)\|_{L^2} \lesssim (1+t)^{-\frac{n}{\sigma}(\frac{1}{m}-\frac{1}{2})+ 1}\|u\|^{p^*}_{X(t)},
\label{t3.6}
\end{equation}
provided that the condition (\ref{spacedimensionCon3}) and the assumption (\ref{existenceCon}) hold. Analogously, we also have
\begin{equation}
\|u_t^{nl}(t,\cdot)\|_{L^2} \lesssim (1+t)^{-\frac{n}{\sigma}(\frac{1}{m}-\frac{1}{2})}\|u\|^{p^*}_{X(t)},
\label{t3.7}
\end{equation}
Now, let us focus on estimating the norm $\big\||D|^r u^{nl}(t,\cdot)\big\|_{L^2}$. Using again the $(L^m \cap L^2)- L^2$ estimates from Proposition \ref{md2.5} we get
\begin{equation}
\big\||D|^r u^{nl}(t,\cdot)\big\|_{L^2} \lesssim \int_0^t (1+t-\tau)^{-\frac{n}{\sigma}(\frac{1}{m}-\frac{1}{2})- \frac{r-\sigma}{\sigma}}\big\||u_t(\tau,\cdot)|^{p^*}\mu\big(|u_t(\tau,\cdot)|\big)\big\|_{L^m \cap L^2 \cap \dot{H}^{r-\sigma}}d\tau.
\label{t3.8}
\end{equation}
Moreover, we can estimate
\begin{equation}
\big\||u_t(\tau,\cdot)|^{p^*}\mu\big(|u_t(\tau,\cdot)|\big)\big\|_{L^m \cap L^2 \cap \dot{H}^{r-\sigma}} \lesssim \big\||u_t(\tau,\cdot)|^{p^*}\big\|_{L^m \cap L^2 \cap \dot{H}^{r-\sigma}} \big\|\mu\big(|u_t(\tau,\cdot)|\big)\big\|_{L^\ity}.
\label{t3.9}
\end{equation}
The estimates for $\big\||u_t(\tau,\cdot)|^{p^*}\big\|_{L^m \cap L^2}$ and $\big\|\mu\big(|u_t(\tau,\cdot)|\big)\big\|_{L^\ity}$ will be handled as in (\ref{t3.3}) and (\ref{t3.4}), respectively. For this reason, we devote our attention to control $\big\||u_t(\tau,\cdot)|^{p^*}\big\|_{\dot{H}^{r-\sigma}}$. Applying Corollary \ref{Corfractionalhomogeneous} for the fractional powers rule and Proposition \ref{ProEmbedding} with $r^* <\frac{n}{2}< r-\sigma$ leads to
\begin{align*}
\big\||u_t(\tau,\cdot)|^{p^*}\big\|_{\dot{H}^{r-\sigma}} &\lesssim \|u_t(\tau,\cdot)\|_{\dot{H}^{r-\sigma}} \|u_t(\tau,\cdot)\|^{p^*-1}_{L^\ity} \\ 
&\lesssim \|u_t(\tau,\cdot)\|_{\dot{H}^{r-\sigma}}\big(\|u_t(\tau,\cdot)\|_{\dot{H}^{r^*}}+ \|u_t(\tau,\cdot)\|_{\dot{H}^{r-\sigma}}\big)^{p^*-1}.
\end{align*}
After employing the fractional Gagliardo-Nirenberg inequality from Proposition \ref{fractionalGagliardoNirenberg}, it deduces
$$ \|u_t(\tau,\cdot)\|_{\dot{H}^{r^*}}\lesssim (1+\tau)^{-\frac{n}{\sigma}(\frac{1}{m}-\frac{1}{2})- \frac{r^*}{\sigma}}\|u\|_{X(\tau)}. $$
Therefore, we arrive at
\begin{equation}
\big\||u_t(\tau,\cdot)|^{p^*}\big\|_{\dot{H}^{r-\sigma}} \lesssim (1+\tau)^{-\frac{np^*}{\sigma}(\frac{1}{m}-\frac{1}{2})- \frac{r-\sigma}{\sigma}- (p^*-1)\frac{r^*}{\sigma}}\|u\|^{p^*}_{X(\tau)} \lesssim (1+\tau)^{-1}\|u\|^{p^*}_{X(\tau)},
\label{t3.10}
\end{equation}
if we choose $r^*= \frac{n}{2}- \e$, where $\e$ is a sufficiently small positive constant. Combining (\ref{t3.3}), (\ref{t3.4}), (\ref{t3.9}) and (\ref{t3.10}) gives
\begin{equation}
\big\||u_t(\tau,\cdot)|^{p^*}\mu\big(|u_t(\tau,\cdot)|\big)\big\|_{L^m \cap L^2 \cap \dot{H}^{r-\sigma}} \lesssim (1+\tau)^{-1} \mu\big(C_1 \e_0 (1+ \tau)^{-\alpha}\big) \|u\|^{p^*}_{X(\tau)}.
\label{t3.11}
\end{equation}
As a result, from (\ref{t3.8}) and (\ref{t3.11}) we get
$$ \big\||D|^r u^{nl}(t,\cdot)\big\|_{L^2}\lesssim  \|u\|^{p^*}_{X(t)}\int_0^t (1+t-\tau)^{-\frac{n}{\sigma}(\frac{1}{m}-\frac{1}{2})- \frac{r-\sigma}{\sigma}}(1+\tau)^{-1} \mu\big(C_1 \e_0 (1+ \tau)^{-\alpha}\big)d\tau. $$
In the analogous way as in the proof of Theorem \ref{dlexistence1} we may conclude the following estimate:
\begin{equation}
\big\||D|^r u^{nl}(t,\cdot)\big\|_{L^2}\lesssim (1+t)^{-\frac{n}{\sigma}(\frac{1}{m}-\frac{1}{2})- \frac{r-\sigma}{\sigma}} \|u\|^{p^*}_{X(t)},
\label{t3.12}
\end{equation}
where the condition (\ref{spacedimensionCon3}) and the assumption (\ref{existenceCon}) are satisfied. Similarly, we also arrive at
\begin{equation}
\big\||D|^{r-\sigma} u_t^{nl}(t,\cdot)\big\|_{L^2}\lesssim (1+t)^{-\frac{n}{\sigma}(\frac{1}{m}-\frac{1}{2})- \frac{r-\sigma}{\sigma}} \|u\|^{p^*}_{X(t)}.
\label{t3.13}
\end{equation}
From all estimates (\ref{t3.6}), (\ref{t3.7}), (\ref{t3.12}), (\ref{t3.13}) and the definition of the norm in $X(t)$, we may conclude immediately the inequality (\ref{pt3.31}). \medskip

Next, let us prove the inequality (\ref{pt3.4}). Then, repeating the proof of Theorem \ref{dlexistence1} and using the same treatment as in the above steps, we may also conclude the inequality (\ref{pt3.4}). Summarizing, the proof is completed.

\subsection{Blow-up results} \label{Sec.Optimal}
The proof of blow-up result in this section is based on a contradiction argument by using the test function method (see, for example, \cite{DabbiccoEbert,Zhang}). In general, this method cannot be directly applied to fractional Laplacian operators $(-\Delta)^\sigma$ and $(-\Delta)^\delta$ as well-known non-local operators. For this reason, the assumption for integers $\sigma$ and $\delta$ comes into play to apply this method. Moreover, due to some techniques in our proofs, we further assume $m= 1$. This means that in this section we only consider the following semi-linear damped $\sigma$-evolution equations:
\begin{equation} \label{pt3.1}
\begin{cases}
u_{tt}+ (-\Delta)^\sigma u+ (-\Delta)^{\delta}u_t= |u|^{p_0}\mu(|u|), \\
u(0,x)= 0,\,\,\, u_t(0,x)=u_1(x),
\end{cases}
\end{equation}
where $\sigma \ge 1$ and $\delta \in [0,\frac{\sigma}{2}]$ are integers, and
\begin{equation} \label{pt3.2}
\begin{cases}
u_{tt}+ (-\Delta)^\sigma u+ (-\Delta)^{\sigma/2}u_t= |u_t|^{p_0}\mu(|u_t|), \\
u(0,x)= 0,\,\,\, u_t(0,x)=u_1(x),
\end{cases}
\end{equation}
with an even integer $\sigma \ge 2$. Here we denote $p_0:= 1+ \frac{2\sigma}{n-2\delta}$ and $p_0:= 1+ \frac{\sigma}{n}$ to (\ref{pt3.1}) and (\ref{pt3.2}), respectively.

\subsubsection{Proof of Theorem \ref{dloptimal1}.}
First, we introduce test function $\varphi= \varphi(r)$ having the following properties:
$$ \varphi \in \mathcal{C}_0^\ity([0,\it)) \text{ and }
\varphi(r)= \begin{cases}
1 \quad \text{ for } r \in [0, 1/2], \\
0 \quad \text{ for } r \in [1,\ity).
\end{cases} $$
Moreover, we assume that $\varphi= \varphi(r)$ is a decreasing function. We also introduce the function $\varphi^*= \varphi^*(r)$ satisfying
$$ \varphi^*(r)= \begin{cases}
0 \qquad\,\, \text{ for } r \in [0, 1/2), \\
\varphi(r) \quad \text{ for } r \in [1/2,\ity).
\end{cases} $$
Let $R$ be a large parameter in $[0,\ity)$. We define the following two functions:
$$ \phi_R(t,x)= \Big(\varphi\Big(\frac{|x|^{2(\sigma-\delta)}+ t}{R}\Big)\Big)^{n+2(\sigma-\delta)}, \text{ and } \phi^*_R(t,x)= \Big(\varphi^*\Big(\frac{|x|^{2(\sigma-\delta)}+ t}{R}\Big)\Big)^{n+2(\sigma-\delta)}. $$
Then it is clear to see that
\begin{align*}
&\text{supp} \phi_R \subset Q_R:= \big\{(t,x): (t,|x|) \in [0,R] \times [0,R^{1/2(\sigma-\delta)}] \big\}, \\ 
&\text{supp} \phi^*_R \subset Q^*_R:= \big\{(t,x): (t,|x|) \in [R/2,R] \times [(R/2)^{1/2\sigma},R^{1/2(\sigma-\delta)}] \big\}.
\end{align*}
Now we define the functional
$$ I_R:= \int_0^{\ity}\int_{\R^n}|u(t,x)|^{p_0}\mu\big(|u(t,x)|\big) \phi_R(t,x)\,dxdt= \int_{Q_R}|u(t,x)|^{p_0}\mu\big(|u(t,x)|\big) \phi_R(t,x)\,d(x,t). $$
Let us assume that $u= u(t,x)$ is a global (in time) Sobolev solution to (\ref{pt3.1}). After multiplying the equation (\ref{pt3.1}) by $\phi_R=\phi_R(t,x)$, we carry out partial integration to derive
\begin{align*}
0\le I_R &= -\int_{\R^n} u_1(x)\phi_R(0,x)\,dx \\
&\qquad+ \int_{Q_R}u(t,x) \big(\partial_t^2 \phi_R(t,x)+ (-\Delta)^{\sigma} \phi_R(t,x)- (-\Delta)^{\delta}\partial_t \phi_R(t,x) \big)\,d(x,t) \\ 
&=: -\int_{\R^n} u_1(x)\phi_R(0,x)\,dx + J_R.
\end{align*}
Because of the assumption (\ref{optimaldata}), there exists a sufficiently large constant $R_0> 0$ such that for all $R > R_0$ it holds
$$ \int_{\R^n} u_1(x)\phi_R(0,x)\,dx >0. $$
Consequently, we obtain
\begin{equation}
0\le I_R < J_R \text{ for all } R > R_0.
\label{t4.1}
\end{equation}
In order to estimate $J_R$, firstly we have
\begin{align}
|\partial_t \phi_R(t,x)| &\lesssim \Big|\frac{1}{R}\Big(\varphi\Big(\frac{|x|^{2(\sigma-\delta)}+ t}{R}\Big)\Big)^{n+2(\sigma-\delta)-1} \varphi'\Big(\frac{|x|^{2(\sigma-\delta)}+ t}{R}\Big)\Big| \nonumber \\ 
&\lesssim \frac{1}{R}\Big(\varphi^*\Big(\frac{|x|^{2(\sigma-\delta)}+ t}{R}\Big)\Big)^{n+2(\sigma-\delta)-1}. \label{t4.2}
\end{align}
Besides, a further calculation leads to
\begin{align}
|\partial_t^2 \phi_R(t,x)| &\lesssim \Big|\frac{1}{R^2}\Big(\varphi\Big(\frac{|x|^{2(\sigma-\delta)}+ t}{R}\Big)\Big)^{n+2(\sigma-\delta)-2} \Big(\varphi'\Big(\frac{|x|^{2(\sigma-\delta)}+ t}{R}\Big)\Big)^2\Big| \nonumber \\
&\qquad + \Big|\frac{1}{R^2}\Big(\varphi\Big(\frac{|x|^{2(\sigma-\delta)}+ t}{R}\Big)\Big)^{n+2(\sigma-\delta)-1} \varphi''\Big(\frac{|x|^{2(\sigma-\delta)}+ t}{R}\Big)\Big| \nonumber \\ 
&\lesssim \frac{1}{R^2}\Big(\varphi^*\Big(\frac{|x|^{2(\sigma-\delta)}+ t}{R}\Big)\Big)^{n+2(\sigma-\delta)-2} \label{t4.3}.
\end{align}
To control $(-\Delta)^{\sigma} \phi_R(t,x)$, we shall apply Lemma \ref{LemmaDerivative} as a main tool. Indeed, we divide our consideration into two sub-steps as follows: \medskip

\noindent Step 1: $\quad$ Applying Lemma \ref{LemmaDerivative} with $h(z)= \varphi(z)$ and $f(x)= \frac{|x|^{2(\sigma-\delta)}+ t}{R}$, we get
\begin{align*}
\Big|\partial_x^\alpha \varphi\Big(\frac{|x|^{2(\sigma-\delta)}+ t}{R}\Big)\Big| &\le \sum_{k=1}^{|\alpha|}\varphi^{(k)} \Big(\frac{|x|^{2(\sigma-\delta)}+ t}{R}\Big) \\
&\qquad \quad \times \Big(\sum_{\substack{|\gamma_1|+\cdots+|\gamma_k|= |\alpha| \\ |\gamma_i|\ge 1}} \Big|\partial_x^{\gamma_1}\Big(\frac{|x|^{2(\sigma-\delta)}+t}{R}\Big)\Big| \cdots \Big|\partial_x^{\gamma_k}\Big(\frac{|x|^{2(\sigma-\delta)}+t}{R}\Big)\Big| \Big) \\
&\le \sum_{k=1}^{|\alpha|}\varphi^{(k)} \Big(\frac{|x|^{2(\sigma-\delta)}+ t}{R}\Big)\Big(\sum_{\substack{|\gamma_1|+\cdots+|\gamma_k|= |\alpha| \\ |\gamma_i|\ge 1}} \frac{|x|^{2(\sigma-\delta)- |\gamma_1|}}{R} \cdots \frac{|x|^{2(\sigma-\delta)- |\gamma_k|}}{R} \Big) \\
&\lesssim \sum_{k=1}^{|\alpha|} \Big(\frac{|x|^{2(\sigma-\delta)}}{R}\Big)^k |x|^{-|\alpha|} \lesssim \frac{|x|^{2(\sigma-\delta)- |\alpha|}}{R}\quad \big(\text{since}\,\, |x|^{2(\sigma-\delta)} \le R \text{ in } Q^*_R\big).
\end{align*}

\noindent Step 2: $\quad$ Applying Lemma \ref{LemmaDerivative} with $h(z)= z^{p_0}$ and $f(x)= \varphi\big(\frac{|x|^{2(\sigma-\delta)}+ t}{R}\big)$, we obtain
\begin{align}
\big|(-\Delta)^{\sigma} \phi_R(t,x)\big| &\le \sum_{k=1}^{2\sigma} \Big(\varphi\Big(\frac{|x|^{2(\sigma-\delta)}+ t}{R}\Big)\Big)^{n+2(\sigma-\delta)- k} \nonumber \\
&\qquad \quad \times \Big(\sum_{\substack{|\gamma_1|+\cdots+|\gamma_k|= 2\sigma \\ |\gamma_i|\ge 1}} \Big|\partial_x^{\gamma_1} \varphi\Big(\frac{|x|^{2(\sigma-\delta)}+ t}{R}\Big)\Big|\, \cdots \, \Big|\partial_x^{\gamma_k} \varphi\Big(\frac{|x|^{2(\sigma-\delta)}+ t}{R} \Big)\Big| \Big) \nonumber \\
&\lesssim \sum_{k=1}^{2\sigma} \Big(\varphi^*\Big(\frac{|x|^{2(\sigma-\delta)}+ t}{R}\Big)\Big)^{n+2(\sigma-\delta)- k} \sum_{\substack{|\gamma_1|+\cdots+|\gamma_k|= 2\sigma \\ |\gamma_i|\ge 1}} \frac{|x|^{2(\sigma-\delta)- |\gamma_1|-\cdots-|\gamma_k|}}{R} \nonumber \\
&\lesssim \sum_{k=1}^{2\sigma} \Big(\varphi^*\Big(\frac{|x|^{2(\sigma-\delta)}+ t}{R}\Big)\Big)^{n+2(\sigma-\delta)- k}\,\, \frac{|x|^{2k(\sigma-\delta)- 2\sigma}}{R^k} \nonumber \\
&\lesssim \frac{1}{R^{\frac{\sigma}{\sigma-\delta}}} \Big(\varphi^*\Big(\frac{|x|^{2(\sigma-\delta)}+ t}{R}\Big)\Big)^{n- 2\delta} \quad \big(\text{since}\,\, R/2\le |x|^{2(\sigma-\delta)} \le R \text{ in } Q^*_R\big). \label{t4.4}
\end{align}
It is clear to see that if $\delta=0$ then $\big|(-\Delta)^{\delta}\partial_t \phi_R(t,x)\big|$ was estimated in (\ref{t4.2}). For the case $\delta \in (0,\frac{\sigma}{2}]$, we can proceed in the analogous way as we controlled $\big|(-\Delta)^{\sigma} \phi_R(t,x)\big|$ to derive
\begin{equation}
\big|(-\Delta)^{\delta}\partial_t \phi_R(t,x)\big| \lesssim \frac{1}{R^{\frac{\sigma}{\sigma-\delta}}} \Big(\varphi^*\Big(\frac{|x|^{2(\sigma-\delta)}+ t}{R}\Big)\Big)^{n+2(\sigma-2\delta)-1}
\label{t4.5}
\end{equation}
From (\ref{t4.2}) to (\ref{t4.5}), we arrive at the following estimate:
\begin{align*}
\big|\partial_t^2 \phi_R(t,x)+ (-\Delta)^{\sigma} \phi_R(t,x)- (-\Delta)^{\delta}\partial_t \phi_R(t,x)\big| &\lesssim \frac{1}{R^{\frac{\sigma}{\sigma-\delta}}} \Big(\varphi^*\Big(\frac{|x|^{2(\sigma-\delta)}+ t}{R}\Big)\Big)^{n-2\delta} \\ 
&= \frac{1}{R^{\frac{\sigma}{\sigma-\delta}}}\big(\phi^*_R(t,x)\big)^{\frac{n-2\delta}{n+2(\sigma-\delta)}}.
\end{align*}
Hence, we may conclude
\begin{equation}
J_R= |J_R| \lesssim \frac{1}{R^{\frac{\sigma}{\sigma-\delta}}} \int_{Q_R}|u(t,x)|\, \big(\phi^*_R(t,x)\big)^{\frac{n-2\delta}{n+2(\sigma-\delta)}}\,d(x,t).
\label{t4.6}
\end{equation}
Now we focus on our attention to estimate the above integral. To do this, we introduce the function $\Psi(s)= s^{p_0}\mu(s)$. Then we derive
\begin{align}
\Psi\Big( |u(t,x)|\, \big(\phi^*_R(t,x)\big)^{\frac{n-2\delta}{n+2(\sigma-\delta)}}\Big) &= |u(t,x)|^{p_0}\, \big(\phi^*_R(t,x)\big)^{\frac{p_0(n-2\delta)}{n+2(\sigma-\delta)}} \mu\Big( |u(t,x)|\, \big(\phi^*_R(t,x)\big)^{\frac{n-2\delta}{n+2(\sigma-\delta)}}\Big) \nonumber \\ 
&\le |u(t,x)|^{p_0}\, \phi^*_R(t,x) \mu\big( |u(t,x)|\big)= \Psi\big( |u(t,x)|\big)\, \phi^*_R(t,x). \label{t4.7}
\end{align}
Here we used the property of the increasing function $\mu(s)$ and the relation
$$0\le \big(\phi^*_R(t,x)\big)^{\frac{n-2\delta}{n+2(\sigma-\delta)}} \le 1. $$
Moreover, it is clear to see that $\Psi$ is a convex function. Applying Proposition \ref{Jenseninequality} with $h(s)= \Psi(s)$, $f(t,x)= |u(t,x)|\big(\phi^*(t,x)\big)^{\frac{n-2\delta}{n+2(\sigma-\delta)}}$ and $\gamma \equiv 1$ gives the following estimate:
$$ \Psi\Big(\frac{\int_{Q^*_R} |u(t,x)|\big(\phi^*(t,x)\big)^{\frac{n-2\delta}{n+2(\sigma-\delta)}}\,d(x,t)}{\int_{Q^*_R} 1\,d(x,t)}\Big) \le \frac{\int_{Q^*_R} \Psi\Big( |u(t,x)|\big(\phi^*_R(t,x)\big)^{\frac{n-2\delta}{n+2(\sigma-\delta)}}\Big)\,d(x,t)}{\int_{Q^*_R} 1\,d(x,t)}. $$
We can easily compute
$$ \int_{Q^*_R} 1\,d(x,t)\approx R^{1+ \frac{n}{2(\sigma-\delta)}}. $$
Hence, we get
\begin{align}
\Psi\Big(\frac{\int_{Q^*_R} |u(t,x)|\big(\phi^*_R(t,x)\big)^{\frac{n-2\delta}{n+2(\sigma-\delta)}}\,d(x,t)}{R^{1+ \frac{n}{2(\sigma-\delta)}}}\Big) &\le \frac{\int_{Q^*_R} \Psi\Big( |u(t,x)|\big(\phi^*(t,x)\big)^{\frac{n-2\delta}{n+2(\sigma-\delta)}}\Big)\,d(x,t)}{R^{1+ \frac{n}{2(\sigma-\delta)}}} \nonumber \\
&\le \frac{\int_{Q_R} \Psi\Big( |u(t,x)|\big(\phi^*_R(t,x)\big)^{\frac{n-2\delta}{n+2(\sigma-\delta)}}\Big)\,d(x,t)}{R^{1+ \frac{n}{2(\sigma-\delta)}}}. \label{t4.8}
\end{align}
Combining the estimates (\ref{t4.7}) and (\ref{t4.8}) we may arrive at
\begin{equation}
\Psi\Big(\frac{\int_{Q^*_R} |u(t,x)|\big(\phi^*_R(t,x)\big)^{\frac{n-2\delta}{n+2(\sigma-\delta)}}\,d(x,t)}{R^{1+ \frac{n}{2(\sigma-\delta)}}}\Big) \le \frac{\int_{Q_R} \Psi\big( |u(t,x)|\big)\, \phi^*_R(t,x)\,d(x,t)}{R^{1+ \frac{n}{2(\sigma-\delta)}}}.
\label{t4.9}
\end{equation}
Since $\mu(s)$ is an increasing function, it immediately follows that $\Psi(s)$ is also increasing function. For this reason, from (\ref{t4.9}) it deduces
\begin{align}
\int_{Q_R} |u(t,x)|\big(\phi^*_R(t,x)\big)^{\frac{n-2\delta}{n+2(\sigma-\delta)}}\,d(x,t)&= \int_{Q^*_R} |u(t,x|\big(\phi^*_R(t,x)\big)^{\frac{n-2\delta}{n+2(\sigma-\delta)}}\,d(x,t) \nonumber \\ 
&\le R^{1+ \frac{n}{2(\sigma-\delta)}}\,\Psi^{-1}\Big(\frac{\int_{Q_R} \Psi\big( |u(t,x)|\big)\, \phi^*_R(t,x)\,d(x,t)}{R^{1+ \frac{n}{2(\sigma-\delta)}}}\Big). \label{t4.10}
\end{align}
From (\ref{t4.1}), (\ref{t4.6}) and (\ref{t4.10}) we may conclude
\begin{equation}
I_R \lesssim R^{\frac{n-2\delta}{2(\sigma-\delta)}}\,\Psi^{-1}\Big(\frac{\int_{Q_R} \Psi\big( |u(t,x)|\big)\, \phi^*_R(t,x)\,d(x,t)}{R^{1+ \frac{n}{2(\sigma-\delta)}}}\Big),
\label{t4.11}
\end{equation}
for all $R> R_0$. Next we introduce the following two functions:
$$ g(r)= \int_{Q_R} \Psi\big( |u(t,x)|\big)\, \phi^*_r(t,x)\,d(x,t), \text{ with } r\in (0,\ity), $$
and
$$ G(R)= \int_0^R g(r)r^{-1}\,dr. $$
Then we re-write
\begin{align*}
G(R) &= \int_0^R \Big(\int_{Q_R} \Psi\big( |u(t,x)|\big)\, \phi^*_r(t,x)\,d(x,t)\Big) r^{-1}\,dr \\
&= \int_{Q_R} \Psi\big( |u(t,x)|\big)\Big(\int_0^R \Big(\varphi^*\Big(\frac{|x|^{2(\sigma-\delta)}+ t}{r}\Big)\Big)^{n+2(\sigma-\delta)} r^{-1}\,dr\Big)\,d(x,t).
\end{align*}
Carrying out change of variables $\tilde{r}= \frac{|x|^{2(\sigma-\delta)}+ t}{r}$ we derive
\begin{align}
G(R)&= \int_{Q_R} \Psi\big( |u(t,x)|\big)\Big(\int_{\frac{|x|^{2(\sigma-\delta)}+ t}{R}}^\ity \big(\varphi^*(\tilde{r})\big)^{n+2(\sigma-\delta)}\, \tilde{r}^{-1}\,d\tilde{r}\Big)\,d(x,t) \nonumber \\ 
&\le \int_{Q_R} \Psi\big( |u(t,x)|\big)\Big(\int_{1/2}^1 \big(\varphi^*(\tilde{r})\big)^{n+2(\sigma-\delta)}\, \tilde{r}^{-1}\,d\tilde{r}\Big)\,d(x,t) \quad \big(\text{since}\,\, \text{supp}\varphi^* \subset [1/2,1]\big) \nonumber \\
&\le \int_{Q_R} \Psi\big( |u(t,x)|\big)\Big(\int_{1/2}^1 \big(\varphi(\tilde{r})\big)^{n+2(\sigma-\delta)}\, \tilde{r}^{-1}\,d\tilde{r}\Big)\,d(x,t) \quad \big(\text{since}\,\, \varphi^* \equiv \varphi \,\,\text{in } [1/2,1]\big) \nonumber \\
&\le \int_{Q_R} \Psi\big( |u(t,x)|\big) \Big(\varphi\Big(\frac{|x|^{2(\sigma-\delta)}+ t}{R}\Big)\Big)^{n+2(\sigma-\delta)}\Big(\int_{1/2}^1 \tilde{r}^{-1}\,d\tilde{r}\Big)\,d(x,t) \quad \big(\text{since}\,\, \varphi \text{ is decreasing}\big) \nonumber \\
&\le \log(1+e) \int_{Q_R} \Psi\big( |u(t,x)|\big) \Big(\varphi\Big(\frac{|x|^{2(\sigma-\delta)}+ t}{R}\Big)\Big)^{n+2(\sigma-\delta)}\,d(x,t)= \log(1+e)I_R. \label{t4.12}
\end{align}
Moreover, it holds the following relation:
\begin{equation}
G'(R)R= g(R)= \int_{Q_R} \Psi\big( |u(t,x)|\big)\, \phi_R^*(t,x)\,d(x,t)
\label{t4.13}
\end{equation}
From (\ref{t4.11}) to (\ref{t4.13}) we get
$$ \frac{G(R)}{\log(1+e)} \le I_R \le C_1 R^{\frac{n-2\delta}{2(\sigma-\delta)}}\,\Psi^{-1}\Big(\frac{G'(R)}{R^{\frac{n}{2(\sigma-\delta)}}}\Big), $$
for all $R> R_0$ and with a suitable positive constant $C_1$. This implies
$$ \Psi\Big(\frac{G(R)}{C_2 R^{\frac{n-2\delta}{2(\sigma-\delta)}}}\Big) \le \frac{G'(R)}{R^{\frac{n}{2(\sigma-\delta)}}}, $$
for all $R> R_0$ and $C_2:= C_1\log(1+e)> 0$. By the definition of the function $\Psi$, the above inequality is equivalent to
$$ \Big(\frac{G(R)}{C_2 R^{\frac{n-2\delta}{2(\sigma-\delta)}}}\Big)^{p_0}\mu\Big(\frac{G(R)}{C_2 R^{\frac{n}{2(\sigma-\delta)}}}\Big) \le \frac{G'(R)}{R^{\frac{n}{2(\sigma-\delta)}}}, $$
for all $R> R_0$. Therefore, we have
$$ \frac{1}{C_3 R}\mu\Big(\frac{G(R)}{C_2 R^{\frac{n-2\delta}{2(\sigma-\delta)}}}\Big) \le \frac{G'(R)}{\big(G(R)\big)^{p_0}}, $$
for all $R> R_0$ and $C_3:= C_2^{p_0}> 0$. Because $G(R)$ is an increasing function, for all $R> R_0$ it immediately follows the following inequality:
$$ \frac{1}{C_3 R}\, \mu\Big(\frac{G(R_0)}{C_2 R^{\frac{n-2\delta}{2(\sigma-\delta)}}}\Big) \le \frac{G'(R)}{\big(G(R)\big)^{p_0}}. $$
By denoting again $\tilde{s}:= R$ and integrating two sides over $[R_0,R]$ we arrive at
\begin{align*}
\frac{1}{C_3} \int_{R_0}^R \frac{1}{\tilde{s}}\, \mu\Big(\frac{1}{C_4 \tilde{s}^{\frac{n-2\delta}{2(\sigma-\delta)}}}\Big)\,d\tilde{s} &\le \int_{R_0}^R\frac{G'(\tilde{s})}{\big(G(\tilde{s})\big)^{p_0}}\,d\tilde{s}= \frac{n-2\delta}{2\sigma} \Big(\frac{1}{\big(G(R_0)\big)^{\frac{2\sigma}{n-2\delta}}}- \frac{1}{\big(G(R)\big)^{\frac{2\sigma}{n-2\delta}}}\Big) \\ 
&\le \frac{n-2\delta}{2\sigma \big(G(R_0)\big)^{\frac{2\sigma}{n-2\delta}}}, 
\end{align*}
where $C_4:= \frac{C_2}{G(R_0)}> 0$. Letting $R \to \ity$ leads to
$$ \frac{1}{C_3} \int_{R_0}^\ity \frac{1}{\tilde{s}}\, \mu\Big(\frac{1}{C_4 \tilde{s}^{\frac{n-2\delta}{2(\sigma-\delta)}}}\Big)\,d\tilde{s} \le \frac{n-2\delta}{2\sigma \big(G(R_0)\big)^{\frac{2\sigma}{n-2\delta}}}. $$
Finally, using change of variables $s= C_4 \tilde{s}^{\frac{n-2\delta}{2(\sigma-\delta)}}$ we may conclude
$$ C\int_{C_0}^\ity \frac{\mu\big(\frac{1}{s}\big)}{s}\,ds \le \frac{n-2\delta}{2\sigma \big(G(R_0)\big)^{\frac{2\sigma}{n-2\delta}}}, $$
where $C:= \frac{2\sigma}{C_3(n-2\delta)}> 0$ and $C_0:= C_4 R_0^{\frac{n-2\delta}{2(\sigma-\delta)}}> 0$ is a sufficiently large constant. This is a contradiction to the assumption (\ref{blowupCon}). Summarizing, the proof is completed.

\subsubsection{Proof of Theorem \ref{dloptimal2}.}
The proof of this theorem is similar to the proof to Theorem \ref{dloptimal1}. For this reason, we only present the steps which are different. Then, we shall repeat some of the arguments as we did in the proof to Theorem \ref{dloptimal1} to conclude the our proof. \medskip

We introduce test function $\varphi$ and $\varphi^*= \varphi^*(r)$ as in Theorem \ref{dloptimal1}. Then, we define two functions:
$$ \phi_R(t,x)= \Big(\varphi\Big(\frac{|x|^{\sigma}+ t}{R}\Big)\Big)^{2(n+\sigma)}, \text{ and } \phi^*_R(t,x)= \Big(\varphi^*\Big(\frac{|x|^{\sigma}+ t}{R}\Big)\Big)^{2(n+\sigma)}. $$
It is clear to see that
\begin{align*}
&\text{supp} \phi_R \subset Q_R:= \big\{(t,x): (t,|x|) \in [0,R] \times [0,R^{1/\sigma}] \big\}, \\ 
&\text{supp} \phi^*_R \subset Q^*_R:= \big\{(t,x): (t,|x|) \in [R/2,R] \times [(R/2)^{1/\sigma},R^{1/\sigma}] \big\}.
\end{align*}
Moreover, we introduce the function
$$ \Phi_R(t,x)= \int_t^\ity \phi_R(\tau,x)d\tau. $$
Because of supp$\phi_R \subset Q_R$, it follows supp$\Phi_R \subset Q_R$. Here we also notice that the relation $\partial_t \Phi_R(t,x)= -\phi_R(t,x)$ holds. Now we define the functional
$$ I_R:= \int_0^{\ity}\int_{\R^n}|u_t(t,x)|^{p_0}\mu\big(|u_t(t,x)|\big) \phi_R(t,x)\,dxdt= \int_{Q_R}|u_t(t,x)|^{p_0}\mu\big(|u_t(t,x)|\big) \phi_R(t,x)\,d(x,t). $$
Let us assume that $u= u(t,x)$ is a global (in time) Sobolev solution to (\ref{pt3.2}). After multiplying the equation (\ref{pt3.2}) by $\phi_R=\phi_R(t,x)$, we carry out partial integration to obtain
\begin{align*}
0\le I_R &= -\int_{\R^n} u_1(x)\phi_R(0,x)\,dx \\
&\qquad + \int_{Q_R}u_t(t,x) \big(-\partial_t \phi_R(t,x)+ (-\Delta)^{\sigma} \Phi_R(t,x)+ (-\Delta)^{\sigma/2} \phi_R(t,x) \big)\,d(x,t) \\ 
&=: -\int_{\R^n} u_1(x)\phi_R(0,x)\,dx + J_R.
\end{align*}
Due to the assumption (\ref{optimaldata2}), there exists a sufficiently large constant $R_0> 0$ such that for all $R > R_0$ it holds
$$ \int_{\R^n} u_1(x)\phi_R(0,x)\,dx >0. $$
As a result, we have
\begin{equation}
0\le I_R < J_R \text{ for all } R > R_0.
\label{t5.1}
\end{equation}
Let us now devote to estimate $J_R$. In the same way as in the proof of Theorem \ref{dloptimal1} we get
$$ |\partial_t \phi_R(t,x)|\lesssim \frac{1}{R}\Big(\varphi^*\Big(\frac{|x|^{\sigma}+ t}{R}\Big)\Big)^{2(n+\sigma)-1}= \frac{1}{R}\big(\phi^*_R(t,x)\big)^{\frac{2(n+\sigma)-1}{2(n+\sigma)}}. $$
In order to control $(-\Delta)^{\sigma} \Phi_R(t,x)$ and $(-\Delta)^{\sigma/2} \phi_R(t,x)$, we shall apply Lemma \ref{LemmaDerivative} as we did in the proof of Theorem \ref{dloptimal1} to conclude the following estimates:
\begin{align*}
\big|(-\Delta)^{\sigma} \phi_R(t,x)\big|& \lesssim \frac{1}{R^2}\Big(\varphi^*\Big(\frac{|x|^{\sigma}+ t}{R}\Big)\Big)^{2n}= \frac{1}{R^2}\big(\phi^*_R(t,x)\big)^{\frac{n}{n+\sigma}}, \\ 
\big|(-\Delta)^{\sigma/2} \phi_R(t,x)\big|& \lesssim \frac{1}{R}\Big(\varphi^*\Big(\frac{|x|^{\sigma}+ t}{R}\Big)\Big)^{2n+\sigma}= \frac{1}{R}\big(\phi^*_R(t,x)\big)^{\frac{2n+\sigma}{2(n+\sigma)}}.
\end{align*}
For this reason, to estimate $(-\Delta)^{\sigma} \Phi_R(t,x)$ we can proceed as follows:
\begin{align*}
\big|(-\Delta)^{\sigma} \Phi_R(t,x)\big|&= \Big|(-\Delta)^{\sigma} \int_t^\ity \phi_R(\tau,x)d\tau\Big|= \Big|\int_t^\ity (-\Delta)^{\sigma} \phi_R(\tau,x)d\tau\Big|\\ 
&\le \int_t^\ity \big|(-\Delta)^{\sigma} \phi_R(\tau,x)\big|d\tau \lesssim \frac{1}{R^2}\int_t^\ity \big(\phi^*_R(\tau,x)\big)^{\frac{n}{n+\sigma}}d\tau= \frac{1}{R^2}\int_t^R \big(\phi^*_R(\tau,x)\big)^{\frac{n}{n+\sigma}}d\tau \\
&\lesssim \frac{1}{R^2}(R-t)\big(\phi^*_R(\tau_0,x)\big)^{\frac{n}{n+\sigma}} \lesssim \frac{1}{R}\big(\phi^*_R(t,x)\big)^{\frac{n}{n+\sigma}}.
\end{align*}
Here we applied the mean value theorem with $\tau_0 \in (t,R)$ and $\phi^*_R$ is decreasing. Therefore, we arrive at
\begin{equation}
\big|-\partial_t \phi_R(t,x)+ (-\Delta)^{\sigma} \Phi_R(t,x)+ (-\Delta)^{\sigma/2} \phi_R(t,x)\big|\lesssim \frac{1}{R}\big(\phi^*_R(t,x)\big)^{\frac{n}{n+\sigma}}.
\label{t5.2}
\end{equation}
From (\ref{t5.1}) and (\ref{t5.2}) we derive for all $R> R_0$
$$ J_R= |J_R| \lesssim \frac{1}{R} \int_{Q_R}|u_t(t,x)|\, \big(\phi^*_R(t,x)\big)^{\frac{n}{n+\sigma}}\,d(x,t). $$
By introducing the functions $\Psi$, $g$ and $G$ as in the proof of Theorem \ref{dloptimal1} and repeating the arguments in this theorem, it deduces the following inequalities:
\begin{equation}
\frac{G(R)}{\log(1+e)} \le I_R \le C_1 R^{\frac{n}{\sigma}}\,\Psi^{-1}\Big(\frac{G'(R)}{R^{\frac{n}{\sigma}}}\Big),
\label{t5.3}
\end{equation}
for all $R> R_0$ and with a suitable positive constant $C_1$. Here we note that
$$ \int_{Q^*_R} 1\,d(x,t)\approx R^{1+ \frac{n}{\sigma}}. $$
By (\ref{t5.3}) a standard calculation as we carried out in the proof of Theorem \ref{dloptimal1} leads to
$$ \frac{1}{C_3 R}\, \mu\Big(\frac{G(R_0)}{C_2 R^{\frac{n}{\sigma}}}\Big) \le \frac{G'(R)}{\big(G(R)\big)^{p_0}}, $$
with some suitable constants $C_2$ and $C_3$. Then, following the analogous treatment as in the proof of Theorem \ref{dloptimal1} gives a contradiction to the assumption (\ref{blowupCon}). Summarizing, the proof is completed. \medskip


\noindent \textbf{Acknowledgments}\medskip

\noindent The PhD study of MSc. T.A. Dao is supported by Vietnamese Government's Scholarship (Grant number: 2015/911).
\medskip

\noindent\textbf{Appendix A}\medskip

\noindent \textit{A.1. Fractional Gagliardo-Nirenberg inequality}

\begin{md} \label{fractionalGagliardoNirenberg}
Let $1<p,\, p_0,\, p_1<\infty$, $\sigma >0$ and $s\in [0,\sigma)$. Then, it holds the following fractional Gagliardo-Nirenberg inequality for all $u\in L^{p_0} \cap \dot{H}^\sigma_{p_1}$:
$$ \|u\|_{\dot{H}^{s}_p}\lesssim \|u\|_{L^{p_0}}^{1-\theta}\,\, \|u\|_{\dot{H}^{\sigma}_{p_1}}^\theta, $$
where $\theta=\theta_{s,\sigma}(p,p_0,p_1)=\frac{\frac{1}{p_0}-\frac{1}{p}+\frac{s}{n}}{\frac{1}{p_0}-\frac{1}{p_1}+\frac{\sigma}{n}}$ and $\frac{s}{\sigma}\leq \theta\leq 1$ .
\end{md}
For the proof one can see \cite{ReissigEbert,Ozawa}.
\medskip

\noindent \textit{A.2. Fractional powers}

\begin{md} \label{PropSickelfractional}
Let $p>1$, $1< r <\infty$ and $u \in H^{s}_r$, where $s \in \big(\frac{n}{r},p\big)$. Let us denote by $F(u)$ one of the functions $|u|^p,\, \pm |u|^{p-1}u$. Then, the following estimate holds:
$$\|F(u)\|_{H^{s}_r}\lesssim \|u\|_{H^{s}_r}\, \|u\|_{L^\infty}^{p-1}.$$
\end{md}

\begin{hq} \label{Corfractionalhomogeneous}
Under the assumptions of Proposition \ref{PropSickelfractional}, it holds: $\| F(u)\|_{\dot{H}^{s}_r}\lesssim \| u\|_{\dot{H}^{s}_r}\, \|u\|_{L^\infty}^{p-1}.$
\end{hq}
The proof can be found in \cite{DuongKainaneReissig,RunSic}.
\medskip

\noindent \textit{A.3. A fractional Sobolev embedding}

\begin{md} \label{ProEmbedding}
Let $0< s_1< \frac{n}{2}< s_2$. Then, for any function $u \in \dot{H}^{s_1} \cap \dot{H}^{s_2}$ we have
\[ \|u\|_{L^\ity} \lesssim \|u\|_{\dot{H}^{s_1}}+ \|u\|_{\dot{H}^{s_2}}. \]
\end{md}
The proof can be found in \cite{DabbiccoEbertLucente}.
\medskip

%

\noindent \textit{A.4. A generalized Jensen's inequality}

\begin{md} \label{Jenseninequality}
Let $\Omega$ be a measurable set respecting a positive measure $\lambda$ so that $\lambda(\Omega)$ is the positive number. Let $f= f(z): \Omega \to \R$ be a $\lambda$-integrable function with the image in $[a,b]$, and let $\gamma= \gamma(z): \Omega \to \R$ be a positive $\lambda$-integrable function. Then each convex function $h= h(s): [a,b] \to \R$ satisfies the following inequality:
$$ h\Big(\frac{\int_\Omega f(z)\gamma(z)\,d\lambda}{\int_\Omega \gamma(z)\,d\lambda}\Big) \le \frac{\int_\Omega h\big(f(z)\big)\gamma(z)\,d\lambda}{\int_\Omega \gamma(z)\,d\lambda}. $$
\end{md}
The proof of this result can be found in \cite{Pavic}.
\medskip

\noindent \textit{A.4. Useful lemma}
\bbd \label{LemmaDerivative}
The following formula of derivative of composed function holds for any multi-index $\alpha$:
$$ \partial_\xi^\alpha h\big(f(\xi)\big)= \sum_{k=1}^{|\alpha|}h^{(k)} \big(f(\xi)\big)\Big(\sum_{\substack{\gamma_1+\cdots+\gamma_k \le \alpha\\ |\gamma_1|+\cdots+|\gamma_k|= |\alpha|,\, |\gamma_i|\ge 1}}\big(\partial_\xi^{\gamma_1} f(\xi)\big) \cdots \big(\partial_\xi^{\gamma_k} f(\xi)\big)\Big), $$
where $h=h(z)$ and $h^{(k)}(z)=\frac{d^k h(z)}{d\,z^k}$.
\ebd
The result can be found in \cite{Simander} at page $202$.
\medskip


\end{document}